\pdfoutput=1
\RequirePackage{snapshot}
\documentclass[]{IEEEtran}
\usepackage[utf8]{inputenc}
\usepackage[T1]{fontenc}

\usepackage{amsmath, amsfonts,amsthm}
\usepackage{bbm}

\usepackage[font=footnotesize,labelfont=bf]{caption}

\usepackage{color}

\usepackage{mathrsfs}

\renewcommand{\hat}{\widehat}
\renewcommand{\tilde}{\widetilde}
\renewcommand{\phi}{\varphi}

\renewcommand{\i}{i\,}

\newcommand{\N}{\mathbb{N}}

\newcommand{\R}{\mathbb{R}}
 
\newcommand{\C}{\mathbb{C}}

\newcommand{\Fc}{\mathcal{F}}

\newcommand{\p}[1]{\left( #1 \right)}

\newcommand{\argmin}{\operatorname{\arg}\min}

\newtheorem{theorem}{Theorem}
\theoremstyle{remark}

\usepackage{aliascnt}
\newaliascnt{conj}{theorem}
\newaliascnt{cor}{theorem}
\newaliascnt{lemma}{theorem}
\newaliascnt{prop}{theorem}
\newaliascnt{definition}{theorem}
\newaliascnt{example}{theorem}
\newaliascnt{notation}{theorem}
\newaliascnt{experiment}{theorem}

\theoremstyle{theorem}

\theoremstyle{definition}

\aliascntresetthe{conj}
\aliascntresetthe{cor}
\aliascntresetthe{lemma}
\aliascntresetthe{prop}
\aliascntresetthe{definition}
\aliascntresetthe{example}
\aliascntresetthe{notation}
\aliascntresetthe{experiment}

\usepackage{snapshot}

\usepackage{booktabs}
\usepackage{multicol}

\usepackage[boxed, ruled, lined]{algorithm2e}

\usepackage[usenames,dvipsnames,svgnames]{xcolor}

\usepackage{enumerate}
\usepackage{csquotes}

\usepackage{subcaption}

\usepackage[pdftex, 
 colorlinks = true,
  linkcolor= black,
  citecolor = black, draft]{hyperref}

\newcommand{\comment}[1]{}

\newcommand{\hl}{}

\usepackage{tikz}
\usepackage{pgfplots}
\usetikzlibrary{arrows,shapes,snakes,automata,backgrounds,petri}

\usepgfplotslibrary{external} 
\pgfkeys{/pgfplots/MyAxisStyle/.style={
        xtick={0, 128, 256},
        font=\footnotesize,
        tick label style = {font=\scriptsize},
        mark size = 1,
        ytick={0, 1},
        xmin=0,
        xmax=256,
        enlargelimits
}}

\newcommand{\dataStyle}{%
\pgfkeys{/pgfplots/MyAxisStyle/.style={
        xtick={0, 128, 256},
        font=\footnotesize,
        tick label style = {font=\scriptsize},
        mark size = 0.7,
        only marks,
        mark=*,
        ytick={0, 0.5, 1},
        xmin=0,
        xmax=256,
        enlargelimits
}}}

\newcommand{\dataStyleIntro}{%
\pgfkeys{/pgfplots/MyAxisStyle/.style={
        xtick={0, 128, 256},
        font=\footnotesize,
        tick label style = {font=\scriptsize},
        mark size = 0.7,
        only marks,
        mark=*,
        ytick={0, 0.5, 1},
        xmin=0,
        xmax=276,
        enlargelimits
}}}

\newcommand{\pottsStyleIntro}{%
\pgfkeys{/pgfplots/MyAxisStyle/.style={
        xtick={0, 128, 256},
        font=\footnotesize,
        tick label style = {font=\scriptsize},
        mark size = 0.7,
        ytick={0, 1},
        xmin=0,
        xmax=276,
        enlargelimits
}}}

\newcommand{\sparsStyle}{%
\pgfkeys{/pgfplots/MyAxisStyle/.style={
        xtick={0, 128, 256},
        font=\footnotesize,
        tick label style = {font=\scriptsize},
        mark size = 1,
        ytick={0, 10},
        xmin=0,
        xmax=256,
        enlargelimits
}}}

\newcommand{\psnrStyle}{%
\pgfkeys{/pgfplots/MyAxisStyle/.style={
        font=\footnotesize,
        tick label style = {font=\scriptsize},
        enlargelimits, 
        	x label style={at={(0.5,0.37)}},
        xlabel = {$\gamma$},
}}}

\newcommand{\logStyle}{%
\pgfkeys{/pgfplots/MyAxisStyle/.style={
        font=\footnotesize,
        tick label style = {font=\scriptsize},
        xlabel = {$\|x\|_0$},
        ylabel = {$\|Ax - b\|_2^2$},
}}}

\pgfkeys{/pgfplots/MyAxisStyleData/.style={
        xtick={0, 128, 256},
        font=\footnotesize,
        tick label style = {font=\scriptsize},
        mark size = 1,
        only marks,
        ytick={0, 0.5, 1},
        xmin=0,
        xmax=256,
        enlargelimits
}}

\title{ 
Jump-sparse and sparse  recovery using 
 Potts functionals
 }
\author{Martin Storath, 
Andreas Weinmann, Laurent Demaret
\thanks{Martin Storath is with the Biomedical Imaging Group, École polytechnique fédérale de Lausanne, Switzerland.}
\thanks{\hl Andreas Weinmann and Laurent Demaret  are both with the Department of Mathematics, Technische Universit\"at M\"unchen, and the Helmholtz Zentrum München, Germany.}
\thanks{
\hl
The research leading to these results has received funding from the European Research Council under the European Union's Seventh Framework Programme (FP7/2007-2013) / ERC grant agreement no.~267439 and the German Federal Ministry for Education and Research under SysTec Grant 0315508.}}
\date{\today}

\begin{document}

\newlength\figureheight
\newlength\figurewidth
\setlength\figureheight{0.15\textwidth}

\maketitle

\begin{abstract}
We recover jump-sparse and sparse signals from blurred
incomplete data  corrupted by (possibly non-Gaussian) noise
using inverse Potts energy functionals.
We obtain analytical results (existence of minimizers, complexity)  on inverse Potts functionals and provide relations to sparsity problems.
We then propose a new optimization method for these functionals which is based on dynamic programming and the alternating direction method of multipliers (ADMM).
A series of experiments shows that the proposed 
method yields very satisfactory jump-sparse and sparse reconstructions, respectively.
We highlight the capability of the method by comparing it with classical and recent approaches such as TV minimization (jump-sparse signals), orthogonal matching pursuit, iterative hard thresholding, and iteratively reweighted $\ell^1$ minimization (sparse signals).
\end{abstract}

\begin{IEEEkeywords}
Jump-sparsity, sparsity, inverse Potts functional, piecewise constant signal,   denoising, deconvolution, incomplete data, segmentation, ADMM. 
\end{IEEEkeywords}

\section{Introduction}

In this article we aim at reconstructing jump-sparse (and sparse) signals $\overline{x} \in \R^n$ 
from linear noisy measurements $b \in \R^m$ (or $\C^m$) given by
\begin{equation*}
	b = A \overline{x} + \text{noise},
\end{equation*}
where $A$ is a (general) $m\times n$ matrix.
The reader may think of $A$ being a Toeplitz matrix modeling blur
or a Fourier matrix, or a combination of both. 
In particular, we deal with  incomplete data meaning that the number of measurements $m$ is significantly smaller than 
the size $n$ of the original signal.
Since this reconstruction problem 
 is in general ill-posed it requires regularization.
This is usually achieved by minimizing a suitable energy functional which expresses a tradeoff between data-fidelity and regularity.
In view of the jump-sparsity of the underlying signal, the number of jumps $\| \nabla  x \|_0   =  |\{ i:  x_i \neq x_{i+1}\}|$ is a natural and powerful regularizing term \cite{potts1952some,mumford1989optimal,boykov2001fast, winkler2002smoothers,candes2008enhancing}.
The corresponding minimization problem, called {\em inverse Potts problem (iPotts),} reads 
\begin{align}\label{eq:inverse_potts}
	P_\gamma(x) = \gamma \, \| \nabla  x \|_0 + \| A x - b \|_p^p \to \mathrm{min}. 
\end{align}
Here the parameter $\gamma > 0$ controls the tradeoff between jump-sparsity and data fidelity
which is measured by some $\ell^p$ norm, $p \geq 1.$
If the noise is Gaussian then $p = 2$ is the natural choice whereas
$p =1$ is the better choice for Laplacian or impulsive noise.
{\hl (We use the notation $F(x) \to \mathrm{min}$ to denote the minimization problem for the functional $F.$)}

The inverse Potts functional is not convex. To avoid the resulting difficulties,
frequently the total variation (TV) penalty $\| \nabla  x\|_1 = \sum_i | x_{i+1} - x_i|$ 
 is used instead for piecewise constant signal restoration
 \cite{rudin1992nonlinear, vogel1996iterative, chambolle1997image, unser2011stochastic,karahanoglu2011signal, needell2013near,needell2013stable}. The TV problem can be solved using convex optimization and the algorithms converge to a global minimum  \cite{chambolle2004algorithm,chambolle2011first,clason2009duality}.
However, the minimizers of the TV problem in general differ from those
of the inverse Potts problem. 
It turns out that minimization of the Potts functional yields genuine jump-sparse signals
whereas TV minimization  does so only approximately, see for instance \autoref{fig:pottsDeconvL1}.

In this work, 
we are first concerned with the question of existence of minimizers
which is more involved than it seems at first glance. In fact, we will see that 
the finite dimensional inverse Potts problem \eqref{eq:inverse_potts} has a minimizer whereas its continuous time counterpart in general need not have a minimizer at all. 
We further show that the inverse Potts problem is NP-hard; 
thus exact minimizers cannot be computed efficiently.
Accepting this fact, we develop an ADMM optimization strategy 
which shows very good recovery performance in practice.
Furthermore, we  shed light on the relation between the jump-sparsity problem \eqref{eq:inverse_potts} and the sparse recovery problem. 
Let us be more precise.

\subsection{Proposed ADMM approach to the inverse Potts problem}

We approach the inverse Potts problem \eqref{eq:inverse_potts} 
using  the \emph{alternating direction method of multipliers (ADMM).} 
ADMM strategies have recently become  very popular in convex optimization especially TV minimization  \cite{wang2008new, ng2010solving, boyd2011distributed, combettes2011proximal,wahlberg2012admm}. They have also shown their usefulness  in non-convex optimization \cite{boyd2011distributed, chartrand2013nonconvex}.
We  propose the iteration 
\begin{equation*}
\left\{
\begin{aligned}
	u^{k+1} &\in \argmin_{u} \gamma \,\| \nabla  u \|_0 + \tfrac{\mu_k}{2} \|  u - (v^k - \tfrac{\lambda^k}{\mu_k}) \|_2^2, 
	\\
	v^{k+1} &= \argmin_{v}  \| Av - b \|_p^p + \tfrac{\mu_k}{2} \| v - (u^{k+1}  + \tfrac{\lambda^k}{\mu_k}) \|_2^2, 
	\\
	\lambda^{k+1} &= \lambda^k + \mu_k(u^{k+1} - v^{k+1}), 
	\end{aligned}
	\right.
\end{equation*}
where the parameter $\mu_k$ is updated by $\mu_{k+1} = \tau\mu_k$ with fixed
 $\tau > 1.$
The key point  is that each subproblem of this \emph{iPotts-ADMM} algorithm is numerically tractable.
The first one is a classical Potts problem (equation \eqref{eq:inverse_potts} with $A = \mathrm{id}$) which can be solved fast and exactly in the univariate case.
For multivariate data, such as images, we use
the strategy of \cite{mumford1989optimal,boykov2001fast}.
The second subproblem consists of minimizing a classical Tikhonov functional.
When $p=2,$ we solve a normal equation and, for $p=1,$ we use a fast semismooth Newton method \cite{clason2010semismooth}.
We further show that our algorithm converges.
Since the  inverse Potts problem is NP-hard, we  cannot expect that it converges to a global minimizer of \eqref{eq:inverse_potts} in general,
but the numerical results are very satisfactory.

\subsection{Inverse Potts problems and sparsity}

The inverse Potts problem is closely connected to the \enquote{Lagrangian formulation} of the sparse recovery problem
\begin{align}\label{eq:sparseLag}
	S_\gamma(x) = \gamma \, \|  x \|_0 + \| A x - b \|_p^p \to \mathrm{min}.
\end{align}
The formulation \eqref{eq:sparseLag} has been considered in \cite{blumensath2009iterative,lai2011unconstrained}, for instance. 
General references concerning sparsity are the books \cite{elad2010sparse,starck2010sparse,mallat2008wavelet} where also a variety of applications may be found in. 

As with the inverse Potts problem and TV minimization, 
one can replace the number  of non-zero entries
$\|  x \|_0$ by the absolute sum $\sum_i|x_i|$
to obtain a convex relaxation of the sparsity problem \eqref{eq:sparseLag} called basis pursuit denoising (BPDN) or $\ell^1$-minimization. 
It is one topic of compressed sensing \cite{baraniuk2007compressive, bruckstein2009sparse,candes2006compressive}
to clarify under which conditions a minimizer of the $\ell^1$-functional
minimizes the sparsity problem \eqref{eq:sparseLag}. 
Positive answers (with a high probability) are obtained under quite restrictive assumptions on the matrix $A$ such as the restricted isometry property \cite{candes2005decoding}. 
If such conditions are not met 
the solutions of BPDN are in general not minimizers of \eqref{eq:sparseLag}. 
Further related work replaces the jump-penalty $\|  x \|_0$ by the non-convex functionals $\|  x \|_q^q$ with $0<q<1$ \cite{lai2011unconstrained, chen2011complexity, chen2012smoothing}.

In this work, instead of using relaxations,  
we transform the sparsity problem \eqref{eq:sparseLag} to an inverse Potts problem of the form \eqref{eq:inverse_potts}.
We show that this can be done for all data fidelity terms based on the $p$-norm with $p \geq 1.$
Thus we may approach the sparsity problem \eqref{eq:sparseLag} using the proposed 
iPotts-ADMM algorithm.

 An approach based on a transformation which is in a certain sense converse to ours is the one in \cite{fornasier2010iterative}. 
There, Blake-Zisserman problems (which are certain discrete Mumford-Shah problems) with $\ell^2$ data terms 
are transformed into separable sparsity type problems which are then approached by iterative thresholding algorithms.

\subsection{Applications and numerical experiments}

We  apply the proposed iPotts-ADMM algorithm to  reconstruct 
jump-sparse signals, which arise in various applications such as stepping rotations of bacterial flagella   \cite{sowa2005direct},
the cross-hybridization of DNA \cite{snijders2001assembly, drobyshev2003specificity, hupe2004analysis}, single-molecule fluorescence resonance energy transfer 
  \cite{joo2008advances}, and
MALDI imaging \cite{schoenmeyer2013automated}.
Here, we recover jump-sparse signals from indirect measurements, for example from blurred data or Fourier data. The measurements are incomplete and corrupted with noise. 
The noise in our examples is Gaussian noise, Laplacian noise, or impulsive noise.
The iPotts-ADMM algorithm is capable of 
recovering jump sparse signals almost perfectly from a reasonable level of noise, and gives in average higher reconstruction qualities than TV minimization.

We further apply the iPotts-ADMM based method to sparse recovery problems,
which for example appear in source localization \cite{malioutov2005sparse} or 
neuroimaging \cite{figueiredo2007gradient}.
As for jump sparse signals we consider blurred data under different types of noise.
In our numerical experiments, we achieve similarly good results as for jump-sparse signals. 
We highlight the capability of our method by comparing it with  
{orthogonal matching pursuit} \cite{petukhov2006fast, temlyakov2003nonlinear,tropp2004greed},
basis pursuit denoising \cite{yang2011alternating}, {iterative hard thresholding} \cite{blumensath2009iterative}
and {iteratively reweighted $\ell^1$ minimization} \cite{candes2008enhancing},
which are the state-of-the-art approaches to sparse recovery.

In order to guarantee reproduciblity an implementation of our algorithms 
is freely available at \url{http://pottslab.de}.

\subsection{Outline of the paper}
We start out to formulate our theoretical results on the inverse Potts problem  in  \autoref{sec:inverse_Potts_relation_sparsity}.
In  \autoref{sec:Potts_ADMM}, we derive an ADMM algorithm for the inverse Potts problem.
In \autoref{sec:experimental_jump_sparse} and \autoref{sec:experimental_sparse}, we provide numerical experiments;
 \autoref{sec:experimental_jump_sparse} deals with jump-sparse signals whereas, in \autoref{sec:experimental_sparse},
we consider sparse signals.
Finally, we supply the proofs in 
\autoref{sec:proofs}.

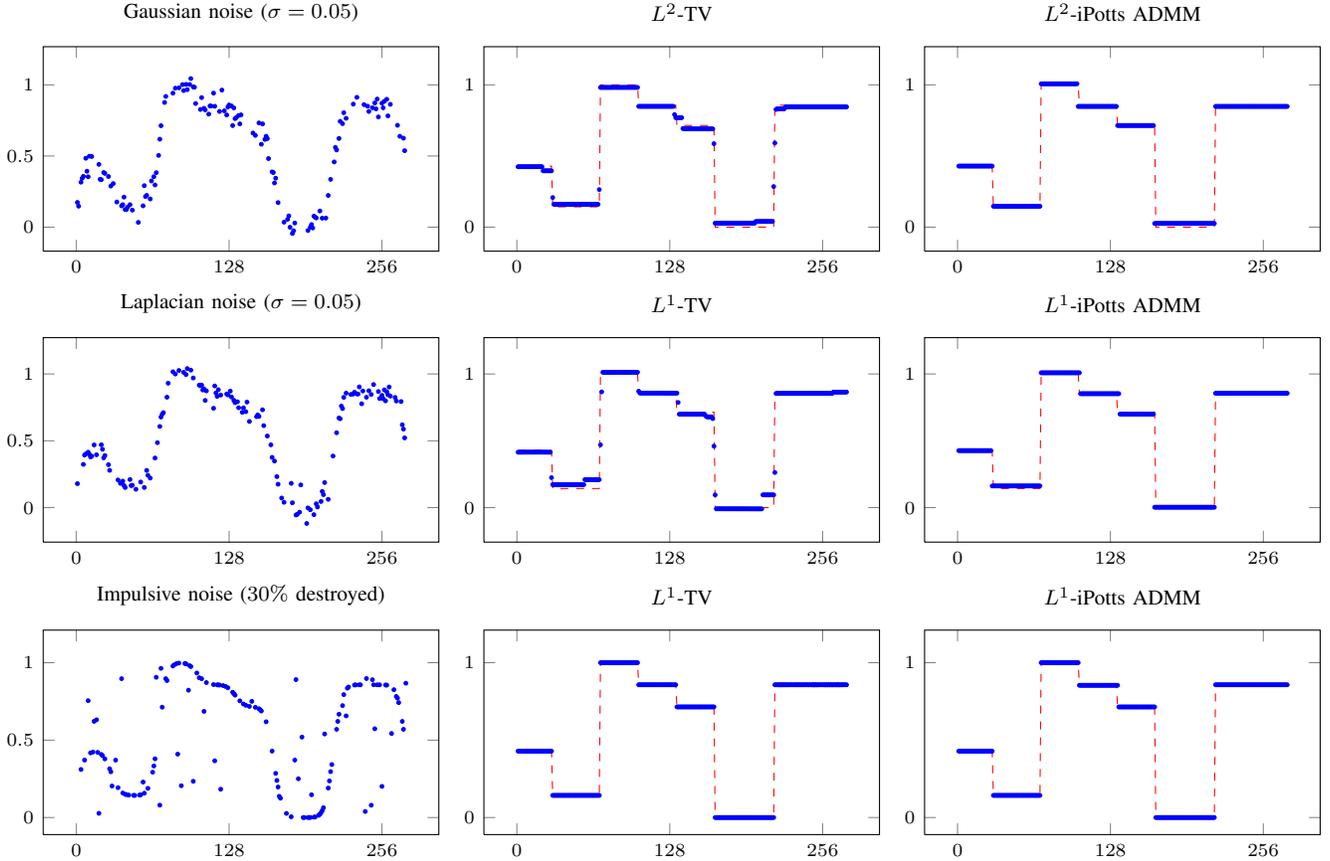
\begin{figure*}
	\def\figfolder{experiments/Potts/}
        \setlength\figurewidth{0.29\textwidth}
        \setlength{\tabcolsep}{-0.5ex}
        \def\rowsep{-0.5ex}
        \begin{tabular}{rrr}
            \\
    		  \def\plottitle{Gaussian noise ($\sigma = 0.05$)}
	\dataStyleIntro
    \def\nextplotname{demoL2PottsInvDeconvData}
    \tikzsetnextfilename{\nextplotname}
%
%
%
\begin{tikzpicture}

\begin{semilogxaxis}[%
view={0}{90},
width=\figurewidth,
height=\figureheight,
scale only axis,
xmin=0.001, xmax=10,
xminorticks=true,
ymin=0, ymax=3,
legend style={at={(0.03,0.97)},anchor=north west,align=left},
MyAxisStyle,
title={\plottitle}]
\addplot [
color=blue,
solid,
mark=o,
mark options={solid}
]
coordinates{
 (0.001,0.100222442857688)(0.00162377673918872,0.130044761898829)(0.00263665089873036,0.163645874121925)(0.0042813323987194,0.191420362715095)(0.00695192796177561,0.224641667058485)(0.0112883789168469,0.263303122036062)(0.0183298071083244,0.317133682911229)(0.0297635144163132,0.397973172790345)(0.0483293023857175,0.527239913146773)(0.0784759970351461,0.711154350279238)(0.127427498570313,0.924204925042849)(0.206913808111479,1.04548772793932)(0.335981828628378,1.18743328122089)(0.545559478116852,1.40628893445207)(0.885866790410082,1.75743527750557)(1.43844988828766,2.45517631726791)(2.33572146909012,2.79991505277133)(3.79269019073225,2.7647742897183)(6.15848211066026,2.83373425054518)(10,2.78463168283138) 
};
\addlegendentry{iPotts-based};

\addplot [
color=green!50!black,
solid,
mark=x,
mark options={solid}
]
coordinates{
 (0.001,0.0896416776129664)(0.00162377673918872,0.130782005099702)(0.00263665089873036,0.171521171292151)(0.0042813323987194,0.201077582134317)(0.00695192796177561,0.237331542179857)(0.0112883789168469,0.285605323916195)(0.0183298071083244,0.357893360816685)(0.0297635144163132,0.469598850019643)(0.0483293023857175,0.646270335731638)(0.0784759970351461,0.849831768663311)(0.127427498570313,0.996468018040267)(0.206913808111479,1.10238250059031)(0.335981828628378,1.30866707331305)(0.545559478116852,1.70102549002036)(0.885866790410082,2.1898689904764)(1.43844988828766,2.79365857081443)(2.33572146909012,2.79991505277133)(3.79269019073225,2.7647742897183)(6.15848211066026,2.83373425054518)(10,2.78463168283138) 
};
\addlegendentry{Direct ADMM};

\end{semilogxaxis}
\end{tikzpicture}%
    &
      \pottsStyleIntro
                \def\plottitle{$L^2$-TV}
    \def\nextplotname{demoL2TVInvDeconvRec}
    \tikzsetnextfilename{\nextplotname}
%
%
%
\begin{tikzpicture}

\begin{semilogxaxis}[%
view={0}{90},
width=\figurewidth,
height=\figureheight,
scale only axis,
xmin=0.001, xmax=10,
xminorticks=true,
ymin=0, ymax=3,
legend style={at={(0.03,0.97)},anchor=north west,align=left},
MyAxisStyle,
title={\plottitle}]
\addplot [
color=blue,
solid,
mark=o,
mark options={solid}
]
coordinates{
 (0.001,0.100222442857688)(0.00162377673918872,0.130044761898829)(0.00263665089873036,0.163645874121925)(0.0042813323987194,0.191420362715095)(0.00695192796177561,0.224641667058485)(0.0112883789168469,0.263303122036062)(0.0183298071083244,0.317133682911229)(0.0297635144163132,0.397973172790345)(0.0483293023857175,0.527239913146773)(0.0784759970351461,0.711154350279238)(0.127427498570313,0.924204925042849)(0.206913808111479,1.04548772793932)(0.335981828628378,1.18743328122089)(0.545559478116852,1.40628893445207)(0.885866790410082,1.75743527750557)(1.43844988828766,2.45517631726791)(2.33572146909012,2.79991505277133)(3.79269019073225,2.7647742897183)(6.15848211066026,2.83373425054518)(10,2.78463168283138) 
};
\addlegendentry{iPotts-based};

\addplot [
color=green!50!black,
solid,
mark=x,
mark options={solid}
]
coordinates{
 (0.001,0.0896416776129664)(0.00162377673918872,0.130782005099702)(0.00263665089873036,0.171521171292151)(0.0042813323987194,0.201077582134317)(0.00695192796177561,0.237331542179857)(0.0112883789168469,0.285605323916195)(0.0183298071083244,0.357893360816685)(0.0297635144163132,0.469598850019643)(0.0483293023857175,0.646270335731638)(0.0784759970351461,0.849831768663311)(0.127427498570313,0.996468018040267)(0.206913808111479,1.10238250059031)(0.335981828628378,1.30866707331305)(0.545559478116852,1.70102549002036)(0.885866790410082,2.1898689904764)(1.43844988828766,2.79365857081443)(2.33572146909012,2.79991505277133)(3.79269019073225,2.7647742897183)(6.15848211066026,2.83373425054518)(10,2.78463168283138) 
};
\addlegendentry{Direct ADMM};

\end{semilogxaxis}
\end{tikzpicture}%
        &
     \pottsStyleIntro
                 \def\plottitle{$L^2$-iPotts ADMM}
    \def\nextplotname{demoL2PottsInvDeconvRec}
    \tikzsetnextfilename{\nextplotname}
%
%
%
\begin{tikzpicture}

\begin{semilogxaxis}[%
view={0}{90},
width=\figurewidth,
height=\figureheight,
scale only axis,
xmin=0.001, xmax=10,
xminorticks=true,
ymin=0, ymax=3,
legend style={at={(0.03,0.97)},anchor=north west,align=left},
MyAxisStyle,
title={\plottitle}]
\addplot [
color=blue,
solid,
mark=o,
mark options={solid}
]
coordinates{
 (0.001,0.100222442857688)(0.00162377673918872,0.130044761898829)(0.00263665089873036,0.163645874121925)(0.0042813323987194,0.191420362715095)(0.00695192796177561,0.224641667058485)(0.0112883789168469,0.263303122036062)(0.0183298071083244,0.317133682911229)(0.0297635144163132,0.397973172790345)(0.0483293023857175,0.527239913146773)(0.0784759970351461,0.711154350279238)(0.127427498570313,0.924204925042849)(0.206913808111479,1.04548772793932)(0.335981828628378,1.18743328122089)(0.545559478116852,1.40628893445207)(0.885866790410082,1.75743527750557)(1.43844988828766,2.45517631726791)(2.33572146909012,2.79991505277133)(3.79269019073225,2.7647742897183)(6.15848211066026,2.83373425054518)(10,2.78463168283138) 
};
\addlegendentry{iPotts-based};

\addplot [
color=green!50!black,
solid,
mark=x,
mark options={solid}
]
coordinates{
 (0.001,0.0896416776129664)(0.00162377673918872,0.130782005099702)(0.00263665089873036,0.171521171292151)(0.0042813323987194,0.201077582134317)(0.00695192796177561,0.237331542179857)(0.0112883789168469,0.285605323916195)(0.0183298071083244,0.357893360816685)(0.0297635144163132,0.469598850019643)(0.0483293023857175,0.646270335731638)(0.0784759970351461,0.849831768663311)(0.127427498570313,0.996468018040267)(0.206913808111479,1.10238250059031)(0.335981828628378,1.30866707331305)(0.545559478116852,1.70102549002036)(0.885866790410082,2.1898689904764)(1.43844988828766,2.79365857081443)(2.33572146909012,2.79991505277133)(3.79269019073225,2.7647742897183)(6.15848211066026,2.83373425054518)(10,2.78463168283138) 
};
\addlegendentry{Direct ADMM};

\end{semilogxaxis}
\end{tikzpicture}%
        \\[\rowsep]
    \dataStyleIntro
            \def\plottitle{Laplacian noise ($\sigma = 0.05$)}
    \def\nextplotname{demoL1PottsInvDeconvData}
    \tikzsetnextfilename{\nextplotname}
%
%
%
\begin{tikzpicture}

\begin{semilogxaxis}[%
view={0}{90},
width=\figurewidth,
height=\figureheight,
scale only axis,
xmin=0.001, xmax=10,
xminorticks=true,
ymin=0, ymax=3,
legend style={at={(0.03,0.97)},anchor=north west,align=left},
MyAxisStyle,
title={\plottitle}]
\addplot [
color=blue,
solid,
mark=o,
mark options={solid}
]
coordinates{
 (0.001,0.100222442857688)(0.00162377673918872,0.130044761898829)(0.00263665089873036,0.163645874121925)(0.0042813323987194,0.191420362715095)(0.00695192796177561,0.224641667058485)(0.0112883789168469,0.263303122036062)(0.0183298071083244,0.317133682911229)(0.0297635144163132,0.397973172790345)(0.0483293023857175,0.527239913146773)(0.0784759970351461,0.711154350279238)(0.127427498570313,0.924204925042849)(0.206913808111479,1.04548772793932)(0.335981828628378,1.18743328122089)(0.545559478116852,1.40628893445207)(0.885866790410082,1.75743527750557)(1.43844988828766,2.45517631726791)(2.33572146909012,2.79991505277133)(3.79269019073225,2.7647742897183)(6.15848211066026,2.83373425054518)(10,2.78463168283138) 
};
\addlegendentry{iPotts-based};

\addplot [
color=green!50!black,
solid,
mark=x,
mark options={solid}
]
coordinates{
 (0.001,0.0896416776129664)(0.00162377673918872,0.130782005099702)(0.00263665089873036,0.171521171292151)(0.0042813323987194,0.201077582134317)(0.00695192796177561,0.237331542179857)(0.0112883789168469,0.285605323916195)(0.0183298071083244,0.357893360816685)(0.0297635144163132,0.469598850019643)(0.0483293023857175,0.646270335731638)(0.0784759970351461,0.849831768663311)(0.127427498570313,0.996468018040267)(0.206913808111479,1.10238250059031)(0.335981828628378,1.30866707331305)(0.545559478116852,1.70102549002036)(0.885866790410082,2.1898689904764)(1.43844988828766,2.79365857081443)(2.33572146909012,2.79991505277133)(3.79269019073225,2.7647742897183)(6.15848211066026,2.83373425054518)(10,2.78463168283138) 
};
\addlegendentry{Direct ADMM};

\end{semilogxaxis}
\end{tikzpicture}%
   &
   \pottsStyleIntro
   \def\plottitle{$L^1$-TV}
   \def\nextplotname{demoL1TVInvDeconvRec}
   \tikzsetnextfilename{\nextplotname}
%
%
%
\begin{tikzpicture}

\begin{semilogxaxis}[%
view={0}{90},
width=\figurewidth,
height=\figureheight,
scale only axis,
xmin=0.001, xmax=10,
xminorticks=true,
ymin=0, ymax=3,
legend style={at={(0.03,0.97)},anchor=north west,align=left},
MyAxisStyle,
title={\plottitle}]
\addplot [
color=blue,
solid,
mark=o,
mark options={solid}
]
coordinates{
 (0.001,0.100222442857688)(0.00162377673918872,0.130044761898829)(0.00263665089873036,0.163645874121925)(0.0042813323987194,0.191420362715095)(0.00695192796177561,0.224641667058485)(0.0112883789168469,0.263303122036062)(0.0183298071083244,0.317133682911229)(0.0297635144163132,0.397973172790345)(0.0483293023857175,0.527239913146773)(0.0784759970351461,0.711154350279238)(0.127427498570313,0.924204925042849)(0.206913808111479,1.04548772793932)(0.335981828628378,1.18743328122089)(0.545559478116852,1.40628893445207)(0.885866790410082,1.75743527750557)(1.43844988828766,2.45517631726791)(2.33572146909012,2.79991505277133)(3.79269019073225,2.7647742897183)(6.15848211066026,2.83373425054518)(10,2.78463168283138) 
};
\addlegendentry{iPotts-based};

\addplot [
color=green!50!black,
solid,
mark=x,
mark options={solid}
]
coordinates{
 (0.001,0.0896416776129664)(0.00162377673918872,0.130782005099702)(0.00263665089873036,0.171521171292151)(0.0042813323987194,0.201077582134317)(0.00695192796177561,0.237331542179857)(0.0112883789168469,0.285605323916195)(0.0183298071083244,0.357893360816685)(0.0297635144163132,0.469598850019643)(0.0483293023857175,0.646270335731638)(0.0784759970351461,0.849831768663311)(0.127427498570313,0.996468018040267)(0.206913808111479,1.10238250059031)(0.335981828628378,1.30866707331305)(0.545559478116852,1.70102549002036)(0.885866790410082,2.1898689904764)(1.43844988828766,2.79365857081443)(2.33572146909012,2.79991505277133)(3.79269019073225,2.7647742897183)(6.15848211066026,2.83373425054518)(10,2.78463168283138) 
};
\addlegendentry{Direct ADMM};

\end{semilogxaxis}
\end{tikzpicture}%
   &
   \pottsStyleIntro
    \def\plottitle{$L^1$-iPotts ADMM}
    \def\nextplotname{demoL1PottsInvDeconvRec}
    \tikzsetnextfilename{\nextplotname}
%
%
%
\begin{tikzpicture}

\begin{semilogxaxis}[%
view={0}{90},
width=\figurewidth,
height=\figureheight,
scale only axis,
xmin=0.001, xmax=10,
xminorticks=true,
ymin=0, ymax=3,
legend style={at={(0.03,0.97)},anchor=north west,align=left},
MyAxisStyle,
title={\plottitle}]
\addplot [
color=blue,
solid,
mark=o,
mark options={solid}
]
coordinates{
 (0.001,0.100222442857688)(0.00162377673918872,0.130044761898829)(0.00263665089873036,0.163645874121925)(0.0042813323987194,0.191420362715095)(0.00695192796177561,0.224641667058485)(0.0112883789168469,0.263303122036062)(0.0183298071083244,0.317133682911229)(0.0297635144163132,0.397973172790345)(0.0483293023857175,0.527239913146773)(0.0784759970351461,0.711154350279238)(0.127427498570313,0.924204925042849)(0.206913808111479,1.04548772793932)(0.335981828628378,1.18743328122089)(0.545559478116852,1.40628893445207)(0.885866790410082,1.75743527750557)(1.43844988828766,2.45517631726791)(2.33572146909012,2.79991505277133)(3.79269019073225,2.7647742897183)(6.15848211066026,2.83373425054518)(10,2.78463168283138) 
};
\addlegendentry{iPotts-based};

\addplot [
color=green!50!black,
solid,
mark=x,
mark options={solid}
]
coordinates{
 (0.001,0.0896416776129664)(0.00162377673918872,0.130782005099702)(0.00263665089873036,0.171521171292151)(0.0042813323987194,0.201077582134317)(0.00695192796177561,0.237331542179857)(0.0112883789168469,0.285605323916195)(0.0183298071083244,0.357893360816685)(0.0297635144163132,0.469598850019643)(0.0483293023857175,0.646270335731638)(0.0784759970351461,0.849831768663311)(0.127427498570313,0.996468018040267)(0.206913808111479,1.10238250059031)(0.335981828628378,1.30866707331305)(0.545559478116852,1.70102549002036)(0.885866790410082,2.1898689904764)(1.43844988828766,2.79365857081443)(2.33572146909012,2.79991505277133)(3.79269019073225,2.7647742897183)(6.15848211066026,2.83373425054518)(10,2.78463168283138) 
};
\addlegendentry{Direct ADMM};

\end{semilogxaxis}
\end{tikzpicture}%
    \\[\rowsep]
    \dataStyleIntro
           \def\plottitle{Impulsive noise ($30 \%$ destroyed)}
    \def\nextplotname{demoL1PottsInvDeconvImpData}
    \tikzsetnextfilename{\nextplotname}
%
%
%
\begin{tikzpicture}

\begin{semilogxaxis}[%
view={0}{90},
width=\figurewidth,
height=\figureheight,
scale only axis,
xmin=0.001, xmax=10,
xminorticks=true,
ymin=0, ymax=3,
legend style={at={(0.03,0.97)},anchor=north west,align=left},
MyAxisStyle,
title={\plottitle}]
\addplot [
color=blue,
solid,
mark=o,
mark options={solid}
]
coordinates{
 (0.001,0.100222442857688)(0.00162377673918872,0.130044761898829)(0.00263665089873036,0.163645874121925)(0.0042813323987194,0.191420362715095)(0.00695192796177561,0.224641667058485)(0.0112883789168469,0.263303122036062)(0.0183298071083244,0.317133682911229)(0.0297635144163132,0.397973172790345)(0.0483293023857175,0.527239913146773)(0.0784759970351461,0.711154350279238)(0.127427498570313,0.924204925042849)(0.206913808111479,1.04548772793932)(0.335981828628378,1.18743328122089)(0.545559478116852,1.40628893445207)(0.885866790410082,1.75743527750557)(1.43844988828766,2.45517631726791)(2.33572146909012,2.79991505277133)(3.79269019073225,2.7647742897183)(6.15848211066026,2.83373425054518)(10,2.78463168283138) 
};
\addlegendentry{iPotts-based};

\addplot [
color=green!50!black,
solid,
mark=x,
mark options={solid}
]
coordinates{
 (0.001,0.0896416776129664)(0.00162377673918872,0.130782005099702)(0.00263665089873036,0.171521171292151)(0.0042813323987194,0.201077582134317)(0.00695192796177561,0.237331542179857)(0.0112883789168469,0.285605323916195)(0.0183298071083244,0.357893360816685)(0.0297635144163132,0.469598850019643)(0.0483293023857175,0.646270335731638)(0.0784759970351461,0.849831768663311)(0.127427498570313,0.996468018040267)(0.206913808111479,1.10238250059031)(0.335981828628378,1.30866707331305)(0.545559478116852,1.70102549002036)(0.885866790410082,2.1898689904764)(1.43844988828766,2.79365857081443)(2.33572146909012,2.79991505277133)(3.79269019073225,2.7647742897183)(6.15848211066026,2.83373425054518)(10,2.78463168283138) 
};
\addlegendentry{Direct ADMM};

\end{semilogxaxis}
\end{tikzpicture}%
    & 
    \pottsStyleIntro
                \def\plottitle{$L^1$-TV}
    \def\nextplotname{demoL1TVInvDeconvImpRec}
    \tikzsetnextfilename{\nextplotname}
%
%
%
\begin{tikzpicture}

\begin{semilogxaxis}[%
view={0}{90},
width=\figurewidth,
height=\figureheight,
scale only axis,
xmin=0.001, xmax=10,
xminorticks=true,
ymin=0, ymax=3,
legend style={at={(0.03,0.97)},anchor=north west,align=left},
MyAxisStyle,
title={\plottitle}]
\addplot [
color=blue,
solid,
mark=o,
mark options={solid}
]
coordinates{
 (0.001,0.100222442857688)(0.00162377673918872,0.130044761898829)(0.00263665089873036,0.163645874121925)(0.0042813323987194,0.191420362715095)(0.00695192796177561,0.224641667058485)(0.0112883789168469,0.263303122036062)(0.0183298071083244,0.317133682911229)(0.0297635144163132,0.397973172790345)(0.0483293023857175,0.527239913146773)(0.0784759970351461,0.711154350279238)(0.127427498570313,0.924204925042849)(0.206913808111479,1.04548772793932)(0.335981828628378,1.18743328122089)(0.545559478116852,1.40628893445207)(0.885866790410082,1.75743527750557)(1.43844988828766,2.45517631726791)(2.33572146909012,2.79991505277133)(3.79269019073225,2.7647742897183)(6.15848211066026,2.83373425054518)(10,2.78463168283138) 
};
\addlegendentry{iPotts-based};

\addplot [
color=green!50!black,
solid,
mark=x,
mark options={solid}
]
coordinates{
 (0.001,0.0896416776129664)(0.00162377673918872,0.130782005099702)(0.00263665089873036,0.171521171292151)(0.0042813323987194,0.201077582134317)(0.00695192796177561,0.237331542179857)(0.0112883789168469,0.285605323916195)(0.0183298071083244,0.357893360816685)(0.0297635144163132,0.469598850019643)(0.0483293023857175,0.646270335731638)(0.0784759970351461,0.849831768663311)(0.127427498570313,0.996468018040267)(0.206913808111479,1.10238250059031)(0.335981828628378,1.30866707331305)(0.545559478116852,1.70102549002036)(0.885866790410082,2.1898689904764)(1.43844988828766,2.79365857081443)(2.33572146909012,2.79991505277133)(3.79269019073225,2.7647742897183)(6.15848211066026,2.83373425054518)(10,2.78463168283138) 
};
\addlegendentry{Direct ADMM};

\end{semilogxaxis}
\end{tikzpicture}%
    & 
    \pottsStyleIntro
                 \def\plottitle{$L^1$-iPotts ADMM}
    \def\nextplotname{demoL1PottsInvDeconvImpRec}
    \tikzsetnextfilename{\nextplotname}
%
%
%
\begin{tikzpicture}

\begin{semilogxaxis}[%
view={0}{90},
width=\figurewidth,
height=\figureheight,
scale only axis,
xmin=0.001, xmax=10,
xminorticks=true,
ymin=0, ymax=3,
legend style={at={(0.03,0.97)},anchor=north west,align=left},
MyAxisStyle,
title={\plottitle}]
\addplot [
color=blue,
solid,
mark=o,
mark options={solid}
]
coordinates{
 (0.001,0.100222442857688)(0.00162377673918872,0.130044761898829)(0.00263665089873036,0.163645874121925)(0.0042813323987194,0.191420362715095)(0.00695192796177561,0.224641667058485)(0.0112883789168469,0.263303122036062)(0.0183298071083244,0.317133682911229)(0.0297635144163132,0.397973172790345)(0.0483293023857175,0.527239913146773)(0.0784759970351461,0.711154350279238)(0.127427498570313,0.924204925042849)(0.206913808111479,1.04548772793932)(0.335981828628378,1.18743328122089)(0.545559478116852,1.40628893445207)(0.885866790410082,1.75743527750557)(1.43844988828766,2.45517631726791)(2.33572146909012,2.79991505277133)(3.79269019073225,2.7647742897183)(6.15848211066026,2.83373425054518)(10,2.78463168283138) 
};
\addlegendentry{iPotts-based};

\addplot [
color=green!50!black,
solid,
mark=x,
mark options={solid}
]
coordinates{
 (0.001,0.0896416776129664)(0.00162377673918872,0.130782005099702)(0.00263665089873036,0.171521171292151)(0.0042813323987194,0.201077582134317)(0.00695192796177561,0.237331542179857)(0.0112883789168469,0.285605323916195)(0.0183298071083244,0.357893360816685)(0.0297635144163132,0.469598850019643)(0.0483293023857175,0.646270335731638)(0.0784759970351461,0.849831768663311)(0.127427498570313,0.996468018040267)(0.206913808111479,1.10238250059031)(0.335981828628378,1.30866707331305)(0.545559478116852,1.70102549002036)(0.885866790410082,2.1898689904764)(1.43844988828766,2.79365857081443)(2.33572146909012,2.79991505277133)(3.79269019073225,2.7647742897183)(6.15848211066026,2.83373425054518)(10,2.78463168283138) 
};
\addlegendentry{Direct ADMM};

\end{semilogxaxis}
\end{tikzpicture}%
        \end{tabular}
    \caption{The original signal  (dashed line) is convolved by a Gaussian kernel of standard deviation $6$ which $\hl m = 138$ measurements are randomly selected from.
The resulting data is corrupted by different types of noise (left).
       The total variation method (TV) mainly reconstructs the constant parts but adds transitional points in between the plateaus  for Gaussian and Laplacian noise.
         The iPotts-ADMM recovers the true signal almost perfectly; in particular, the correct number of jumps. 
   For impulsive noise, the iPotts-ADMM and the TV method perform equally well.
    }
    \label{fig:pottsDeconvL1}
\end{figure*}

\section{Inverse Potts problems and their relation to sparsity}
\label{sec:inverse_Potts_relation_sparsity}

We start our analysis of the inverse Potts problem by considering the question of existence of minimizers. 
It is remarkable that there is a significant difference between the finite dimensional discrete time case and its
infinite dimensional continuous time counterpart. More precisely, 
we obtain a positive answer for the discrete time problem \eqref{eq:inverse_potts}
but a negative answer for the corresponding continuous time problem.
\begin{theorem}\label{thm:minimizer_existence}
	The inverse Potts problem \eqref{eq:inverse_potts} has a minimizer.
\end{theorem}
The proof of \autoref{thm:minimizer_existence} is given in \autoref{ssec:existMin}.
It uses the compactness of the closed unit ball and the lower boundedness
of an injective linear mapping which are features of finite dimensional spaces. 
Thus it does not carry over to the infinite dimensional continuous time case. 
{\hl
We note that the existence of minimizers for Blake-Zisserman functionals with $\ell^2$ data term has been shown in \cite{fornasier2010iterative}.
For $\ell^2$ data terms, modifications of the proofs of \cite{fornasier2010iterative} would also apply to our setting.
However, for general $\ell^p$ data term, the approach of \cite{fornasier2010iterative} does not carry over.}

The next theorem states that the continuous time counterpart of \autoref{thm:minimizer_existence} is false in general. 
The continuous time counterpart of \eqref{eq:inverse_potts} is obtained by replacing  
the finite dimensional signal and data spaces by $L^p$ function spaces and the matrix $A$
by a bounded operator $A$ between those function spaces.
\begin{theorem}\label{thm:minimizer_not_existence_continuous_time}
	There are linear operators $A$ and data $b$ in $L^p,$ $\hl 1 \leq p < \infty,$  such that the continuous time inverse Potts problem with respect to $A$ and $b$ does not have a minimizer.
\end{theorem}
The proof of \autoref{thm:minimizer_not_existence_continuous_time} is given in \autoref{ssec:existMin}.
The explicit counter-examples we give are convolution operators
 which are in fact important from a practical point of view.  

The next natural step after showing the existence of minimizers (in the discrete case) is to clarify the 
complexity of computing such a minimizer.
We obtain the following result.
\begin{theorem}\label{thm:Potts_NP_hard}
	The inverse Potts problem \eqref{eq:inverse_potts} is NP hard.
\end{theorem}
As a consequence, a fast exact algorithm is not available (unless $P = NP$)
and one has to resort to approximative strategies (see \autoref{sec:Potts_ADMM}). 
The proof of \autoref{thm:Potts_NP_hard} is given in \autoref{ssec:RelSparse}.

Finally, we are interested in the relations between sparsity problems and {\hl univariate} inverse Potts problems.
We first consider the sparsity problem \eqref{eq:sparseLag}. We find a corresponding 
{\hl univariate} inverse Potts problem  whose minimizers are directly related 
to the minimizers of the initial sparsity problem.  
{\hl We use this relation in Section V to apply our algorithm to sparsity problems.}
\begin{theorem}
    \label{thm:Pottsforsparse} Let 
$x^* \in\R^{n+1} $ be a minimizer of the inverse Potts functional associated with the matrix $B =A \nabla$, i.e.,
\begin{equation} \label{eq:inverse_potts2}
  x^* \in \argmin_{x\in \R^{n+1}}  \gamma \|\nabla x\|_0 + \|B x-b \|^p_p. 
\end{equation}
Then $u^* = \nabla x^*$ minimizes the sparsity problem \eqref{eq:sparseLag} related to the matrix $A$ and data $b$.
\end{theorem}

We obtain a converse result for $p=2$ {\hl(still for the univariate setting)}.
The relations between the matrices $A$ and $B$ and between the data
can be given explicitly but are not as simple as above.
{\hl A similar relation has been used in \cite{fornasier2010iterative}
in the context of Blake-Zisserman functionals.} 
The construction does not work for general $p \neq 2$
and it is not clear to us how to get a converse result when $p \neq 2.$ 

\begin{theorem}
   \label{thm:sparseforPotts}
For the inverse Potts problem \eqref{eq:inverse_potts} associated with the matrix $A$ and data $b$
we consider the sparsity problem associated with the matrix $B = A^\prime \nabla^+$ and data $b^\prime.$ {\hl Here $\nabla^+$ is the pseudo-inverse of the discrete difference operator given by \eqref{eq:defNablaPlus}. The modified data $A'$ and $b'$ are given in terms of $A$ and $b$ by \eqref{eq:Aprime} and \eqref{eq:bprime}, respectively.}
Let $u^*$ be a minimizer of the sparsity problem with respect to $B,b'$,i. e., 
\begin{equation} \label{eq:sparse2}
   u^* \in \argmin_{u \in \R^{n-1}}  \gamma \|u\|_0 + \|B u - b^\prime\|_2^2.
\end{equation}
Then $x^* = \nabla^+ u^* + \mu (\nabla^+ u^*) e$ (with $\mu$ given by \eqref{eq:mu_of_x0}) is a solution of the inverse Potts problem \eqref{eq:inverse_potts} associated with $A,b$. 
\end{theorem}

The proofs of \autoref{thm:Pottsforsparse} and \autoref{thm:sparseforPotts} are given in \autoref{ssec:RelSparse}.

\section{Minimization of the Potts functional using the alternating direction method of multipliers}
\label{sec:Potts_ADMM}

In this section, we present our iterative approach to the inverse Potts problem \eqref{eq:inverse_potts}.

\subsection{A new ADMM algorithm for the inverse Potts problem}

The inverse Potts problem is equivalent to the bivariate constrained optimization problem
\begin{equation}\label{eq:potts_split}
\begin{split}
	&\text{minimize } \quad \gamma \| \nabla  u \|_0 + \| Av - b \|_p^p  \\
	&\text{subject to } \quad u - v = 0.
	\end{split}
\end{equation}
We incorporate the constraint $u-v$ into the target functional to obtain 
the unconstrained problem 
 \begin{equation}\label{eq:lagrangian}
\begin{split}
 	L_\mu(u,v,\lambda) = \, &  \gamma \|\nabla  u \|_0 + \langle \lambda , u - v \rangle  \\
	& +  \tfrac{\mu}{2}  \| u-v  \|_2^2 + \| A v - b \|_p^p \to \mathrm{min.}
	\end{split}
\end{equation}
The parameter $\mu > 0$ regulates the coupling of $u$ and $v.$
The dual variable $\lambda$ is an $n$-dimensional vector of Lagrange multipliers.
Equation \eqref{eq:lagrangian} is called the \emph{augmented Lagrangian} of \eqref{eq:potts_split}.
Completing the square in the second and third term of \eqref{eq:lagrangian} yields
\begin{equation}\label{eq:lagrangian_se}
\begin{split}
	L_\mu(u,v,\lambda) = &  \gamma \|\nabla  u \|_0 - \tfrac{\mu}{2} \| \tfrac\lambda\mu \|_2^2  \\
	& +  \tfrac{\mu}{2} \| u-v + \tfrac\lambda\mu \|_2^2 + \| Av - b \|_p^p.
	\end{split}
\end{equation}
In order to minimize the augmented Lagrangian \eqref{eq:lagrangian_se}
we use the \emph{alternating direction method of multipliers (ADMM)}, see e.g.  \cite{boyd2011distributed}.
In the ADMM iteration we first fix $v$ and $\lambda$ and minimize $L_\mu(u,v,\lambda)$  with respect to $u.$ Then we  minimize $L_\mu(u,v,\lambda)$ with respect to $v,$ keeping $u$ and $\lambda$ fixed.
The third step is the update of the dual variable $\lambda.$ 
Thus, the alternating direction method of multipliers for the inverse Potts problem \eqref{eq:inverse_potts} reads
\begin{equation}
	\label{eq:ADMM}
\left\{
\begin{aligned}
	u^{k+1} &\in \argmin_{u} \gamma \| \nabla  u \|_0 + \tfrac{\mu}{2} \|  u - (v^k - \tfrac{\lambda^k}{\mu}) \|_2^2, 
	\\
	v^{k+1} &= \argmin_{v}  \| Av - b \|_p^p + \tfrac{\mu}{2} \| v - (u^{k+1}  + \tfrac{\lambda^k}{\mu}) \|_2^2, 
	\\
	\lambda^{k+1} &= \lambda^k + \mu (u^{k+1} - v^{k+1}).
	\end{aligned}
	\right.
\end{equation}
The crucial point is that both subproblems appearing in the first and the second line of \eqref{eq:ADMM} are computationally tractable {\hl (for $p \in [1, \infty]$)}. The first subproblem is the  minimization of a classical Potts problem
which we elaborate on in \autoref{ssec:classical_potts}.
The second subproblem is the minimization of a classical Tikhonov-type problem
which we explain in \autoref{ssec:tikhonov}.

We initialize the iteration with a small positive coupling parameter $\mu_0 > 0$
and increase it during the iteration by a factor $\tau > 1.$
Hence, $\mu$ is given by the geometric progression
\[
	\mu = \mu_k = \tau^k \cdot \mu_0.
\]
This assures that $u$ and $v$ can evolve quite independently at the beginning and that they are close to each other at the end of the iteration.
We stop the iteration when the norm of $u - v$ falls below some tolerance.
Our approach to the inverse Potts problem is summed up in Algorithm \ref{alg:admm}.
\begin{algorithm}[ht]
\small
	 \KwIn{Data $b \in \R^m,$ model parameter $\gamma > 0,$ measurement matrix $A \in \R^{m \times n}$ }
	 \KwOut{Computed result $u \in \R^n$ to the inverse Potts problem \eqref{eq:inverse_potts} }
	 \Begin{
	 $v \leftarrow {\hl A^* b}$;
	 $\mu \leftarrow \mu_0$;
	 $\lambda \leftarrow 0$\;
 	\Repeat{$\| u - v \|_2^2  < \rm{TOL}$}{
	 \SetKwHangingKw{Let}{$\leftarrow$}
	 $u$ \Let{
	   Minimizer of classical $L^2$-Potts functional \eqref{eq:non_inverse_potts} with data $d=v - \frac\lambda\mu$ and parameter $\delta=\frac{ 2\gamma}{\mu}$}	
	   	$v$ \Let{ Solution 
		of  Tikhonov problem \eqref{eq:tikhonov} 
		with data $b,$ offset vector $w = u + \frac{\lambda}{\mu}$ and parameter $\tfrac{\mu}2$}
		$\lambda$ \Let{$ \lambda + \mu(u - v)$}
		$\mu$ $\leftarrow \tau\cdot \mu$\;
	}
	}
	\caption{iPotts-ADMM}
	\label{alg:admm}
 \end{algorithm}
 
We have the following convergence result, whose proof is given in \autoref{ssec:ProofConv}.
 \begin{theorem}\label{thm:convergence} 
 The ADMM iteration \eqref{eq:ADMM}, and thus Algorithm \ref{alg:admm}, converges.
\end{theorem}
Although we cannot expect convergence to a global minimum for the NP-hard inverse Potts problem, we see in the experimental section that Algorithm \ref{alg:admm} gives very satisfactory reconstruction results.

In our experiments, reasonable numerical values for the parameters in Algorithm \ref{alg:admm} are
 $\mu_0 = \gamma \cdot 10^{-6}$ as initial coupling,
$\tau = 1.05$  for the increment of the coupling,
 and $\mathrm{TOL} = 10^{-6}$  for the stopping tolerance.

\subsection{Minimization of the classical Potts subproblem}
\label{ssec:classical_potts}

 \begin{figure*}
\hspace*{-5ex}        
\setlength\figurewidth{0.15\textwidth}
\setlength{\tabcolsep}{-0.5ex}
\def\folder{experiments/Potts/fourierPSNRPlots/}
        \begin{tabular}{rrrrr}
            \def\plottitle{Data (absolute value)}
    \def\nextplotname{demoL2PottsInvFourierData}
\pgfkeys{/pgfplots/MyAxisStyle/.style={
        xtick={0, 128, 256},
        font=\footnotesize,
        tick label style = {font=\scriptsize},
        mark size = 1,
        only marks,
        ytick={1, 10},
        xmin=0,
        xmax=256,
        enlargelimits
}}
   \tikzsetnextfilename{\nextplotname}
%
%
%
\begin{tikzpicture}

\begin{semilogxaxis}[%
view={0}{90},
width=\figurewidth,
height=\figureheight,
scale only axis,
xmin=0.001, xmax=10,
xminorticks=true,
ymin=0, ymax=3,
legend style={at={(0.03,0.97)},anchor=north west,align=left},
MyAxisStyle,
title={\plottitle}]
\addplot [
color=blue,
solid,
mark=o,
mark options={solid}
]
coordinates{
 (0.001,0.100222442857688)(0.00162377673918872,0.130044761898829)(0.00263665089873036,0.163645874121925)(0.0042813323987194,0.191420362715095)(0.00695192796177561,0.224641667058485)(0.0112883789168469,0.263303122036062)(0.0183298071083244,0.317133682911229)(0.0297635144163132,0.397973172790345)(0.0483293023857175,0.527239913146773)(0.0784759970351461,0.711154350279238)(0.127427498570313,0.924204925042849)(0.206913808111479,1.04548772793932)(0.335981828628378,1.18743328122089)(0.545559478116852,1.40628893445207)(0.885866790410082,1.75743527750557)(1.43844988828766,2.45517631726791)(2.33572146909012,2.79991505277133)(3.79269019073225,2.7647742897183)(6.15848211066026,2.83373425054518)(10,2.78463168283138) 
};
\addlegendentry{iPotts-based};

\addplot [
color=green!50!black,
solid,
mark=x,
mark options={solid}
]
coordinates{
 (0.001,0.0896416776129664)(0.00162377673918872,0.130782005099702)(0.00263665089873036,0.171521171292151)(0.0042813323987194,0.201077582134317)(0.00695192796177561,0.237331542179857)(0.0112883789168469,0.285605323916195)(0.0183298071083244,0.357893360816685)(0.0297635144163132,0.469598850019643)(0.0483293023857175,0.646270335731638)(0.0784759970351461,0.849831768663311)(0.127427498570313,0.996468018040267)(0.206913808111479,1.10238250059031)(0.335981828628378,1.30866707331305)(0.545559478116852,1.70102549002036)(0.885866790410082,2.1898689904764)(1.43844988828766,2.79365857081443)(2.33572146909012,2.79991505277133)(3.79269019073225,2.7647742897183)(6.15848211066026,2.83373425054518)(10,2.78463168283138) 
};
\addlegendentry{Direct ADMM};

\end{semilogxaxis}
\end{tikzpicture}%
    &
         \psnrStyle
        \def\plottitle{PSNR of TV solution}
        \def\nextplotname{demoL2TVInvFourierPSNR}
           \tikzsetnextfilename{\nextplotname}
%
%
%
\begin{tikzpicture}

\begin{semilogxaxis}[%
view={0}{90},
width=\figurewidth,
height=\figureheight,
scale only axis,
xmin=0.001, xmax=10,
xminorticks=true,
ymin=0, ymax=3,
legend style={at={(0.03,0.97)},anchor=north west,align=left},
MyAxisStyle,
title={\plottitle}]
\addplot [
color=blue,
solid,
mark=o,
mark options={solid}
]
coordinates{
 (0.001,0.100222442857688)(0.00162377673918872,0.130044761898829)(0.00263665089873036,0.163645874121925)(0.0042813323987194,0.191420362715095)(0.00695192796177561,0.224641667058485)(0.0112883789168469,0.263303122036062)(0.0183298071083244,0.317133682911229)(0.0297635144163132,0.397973172790345)(0.0483293023857175,0.527239913146773)(0.0784759970351461,0.711154350279238)(0.127427498570313,0.924204925042849)(0.206913808111479,1.04548772793932)(0.335981828628378,1.18743328122089)(0.545559478116852,1.40628893445207)(0.885866790410082,1.75743527750557)(1.43844988828766,2.45517631726791)(2.33572146909012,2.79991505277133)(3.79269019073225,2.7647742897183)(6.15848211066026,2.83373425054518)(10,2.78463168283138) 
};
\addlegendentry{iPotts-based};

\addplot [
color=green!50!black,
solid,
mark=x,
mark options={solid}
]
coordinates{
 (0.001,0.0896416776129664)(0.00162377673918872,0.130782005099702)(0.00263665089873036,0.171521171292151)(0.0042813323987194,0.201077582134317)(0.00695192796177561,0.237331542179857)(0.0112883789168469,0.285605323916195)(0.0183298071083244,0.357893360816685)(0.0297635144163132,0.469598850019643)(0.0483293023857175,0.646270335731638)(0.0784759970351461,0.849831768663311)(0.127427498570313,0.996468018040267)(0.206913808111479,1.10238250059031)(0.335981828628378,1.30866707331305)(0.545559478116852,1.70102549002036)(0.885866790410082,2.1898689904764)(1.43844988828766,2.79365857081443)(2.33572146909012,2.79991505277133)(3.79269019073225,2.7647742897183)(6.15848211066026,2.83373425054518)(10,2.78463168283138) 
};
\addlegendentry{Direct ADMM};

\end{semilogxaxis}
\end{tikzpicture}%
    &
             \psnrStyle
        \def\plottitle{PSNR of iPotts solution}
        \def\nextplotname{demoL2PottsInvFourierPSNR}
           \tikzsetnextfilename{\nextplotname}
%
%
%
\begin{tikzpicture}

\begin{semilogxaxis}[%
view={0}{90},
width=\figurewidth,
height=\figureheight,
scale only axis,
xmin=0.001, xmax=10,
xminorticks=true,
ymin=0, ymax=3,
legend style={at={(0.03,0.97)},anchor=north west,align=left},
MyAxisStyle,
title={\plottitle}]
\addplot [
color=blue,
solid,
mark=o,
mark options={solid}
]
coordinates{
 (0.001,0.100222442857688)(0.00162377673918872,0.130044761898829)(0.00263665089873036,0.163645874121925)(0.0042813323987194,0.191420362715095)(0.00695192796177561,0.224641667058485)(0.0112883789168469,0.263303122036062)(0.0183298071083244,0.317133682911229)(0.0297635144163132,0.397973172790345)(0.0483293023857175,0.527239913146773)(0.0784759970351461,0.711154350279238)(0.127427498570313,0.924204925042849)(0.206913808111479,1.04548772793932)(0.335981828628378,1.18743328122089)(0.545559478116852,1.40628893445207)(0.885866790410082,1.75743527750557)(1.43844988828766,2.45517631726791)(2.33572146909012,2.79991505277133)(3.79269019073225,2.7647742897183)(6.15848211066026,2.83373425054518)(10,2.78463168283138) 
};
\addlegendentry{iPotts-based};

\addplot [
color=green!50!black,
solid,
mark=x,
mark options={solid}
]
coordinates{
 (0.001,0.0896416776129664)(0.00162377673918872,0.130782005099702)(0.00263665089873036,0.171521171292151)(0.0042813323987194,0.201077582134317)(0.00695192796177561,0.237331542179857)(0.0112883789168469,0.285605323916195)(0.0183298071083244,0.357893360816685)(0.0297635144163132,0.469598850019643)(0.0483293023857175,0.646270335731638)(0.0784759970351461,0.849831768663311)(0.127427498570313,0.996468018040267)(0.206913808111479,1.10238250059031)(0.335981828628378,1.30866707331305)(0.545559478116852,1.70102549002036)(0.885866790410082,2.1898689904764)(1.43844988828766,2.79365857081443)(2.33572146909012,2.79991505277133)(3.79269019073225,2.7647742897183)(6.15848211066026,2.83373425054518)(10,2.78463168283138) 
};
\addlegendentry{Direct ADMM};

\end{semilogxaxis}
\end{tikzpicture}%
        &
             \def\plottitle{$L^2$-TV}
    \def\nextplotname{demoL2TVInvFourierRec}
    \tikzsetnextfilename{\nextplotname}
%
%
%
\begin{tikzpicture}

\begin{semilogxaxis}[%
view={0}{90},
width=\figurewidth,
height=\figureheight,
scale only axis,
xmin=0.001, xmax=10,
xminorticks=true,
ymin=0, ymax=3,
legend style={at={(0.03,0.97)},anchor=north west,align=left},
MyAxisStyle,
title={\plottitle}]
\addplot [
color=blue,
solid,
mark=o,
mark options={solid}
]
coordinates{
 (0.001,0.100222442857688)(0.00162377673918872,0.130044761898829)(0.00263665089873036,0.163645874121925)(0.0042813323987194,0.191420362715095)(0.00695192796177561,0.224641667058485)(0.0112883789168469,0.263303122036062)(0.0183298071083244,0.317133682911229)(0.0297635144163132,0.397973172790345)(0.0483293023857175,0.527239913146773)(0.0784759970351461,0.711154350279238)(0.127427498570313,0.924204925042849)(0.206913808111479,1.04548772793932)(0.335981828628378,1.18743328122089)(0.545559478116852,1.40628893445207)(0.885866790410082,1.75743527750557)(1.43844988828766,2.45517631726791)(2.33572146909012,2.79991505277133)(3.79269019073225,2.7647742897183)(6.15848211066026,2.83373425054518)(10,2.78463168283138) 
};
\addlegendentry{iPotts-based};

\addplot [
color=green!50!black,
solid,
mark=x,
mark options={solid}
]
coordinates{
 (0.001,0.0896416776129664)(0.00162377673918872,0.130782005099702)(0.00263665089873036,0.171521171292151)(0.0042813323987194,0.201077582134317)(0.00695192796177561,0.237331542179857)(0.0112883789168469,0.285605323916195)(0.0183298071083244,0.357893360816685)(0.0297635144163132,0.469598850019643)(0.0483293023857175,0.646270335731638)(0.0784759970351461,0.849831768663311)(0.127427498570313,0.996468018040267)(0.206913808111479,1.10238250059031)(0.335981828628378,1.30866707331305)(0.545559478116852,1.70102549002036)(0.885866790410082,2.1898689904764)(1.43844988828766,2.79365857081443)(2.33572146909012,2.79991505277133)(3.79269019073225,2.7647742897183)(6.15848211066026,2.83373425054518)(10,2.78463168283138) 
};
\addlegendentry{Direct ADMM};

\end{semilogxaxis}
\end{tikzpicture}%
    &
                  \def\plottitle{$L^2$-iPotts ADMM}
    \def\nextplotname{demoL2PottsInvFourierRec}
    \tikzsetnextfilename{\nextplotname}
%
%
%
\begin{tikzpicture}

\begin{semilogxaxis}[%
view={0}{90},
width=\figurewidth,
height=\figureheight,
scale only axis,
xmin=0.001, xmax=10,
xminorticks=true,
ymin=0, ymax=3,
legend style={at={(0.03,0.97)},anchor=north west,align=left},
MyAxisStyle,
title={\plottitle}]
\addplot [
color=blue,
solid,
mark=o,
mark options={solid}
]
coordinates{
 (0.001,0.100222442857688)(0.00162377673918872,0.130044761898829)(0.00263665089873036,0.163645874121925)(0.0042813323987194,0.191420362715095)(0.00695192796177561,0.224641667058485)(0.0112883789168469,0.263303122036062)(0.0183298071083244,0.317133682911229)(0.0297635144163132,0.397973172790345)(0.0483293023857175,0.527239913146773)(0.0784759970351461,0.711154350279238)(0.127427498570313,0.924204925042849)(0.206913808111479,1.04548772793932)(0.335981828628378,1.18743328122089)(0.545559478116852,1.40628893445207)(0.885866790410082,1.75743527750557)(1.43844988828766,2.45517631726791)(2.33572146909012,2.79991505277133)(3.79269019073225,2.7647742897183)(6.15848211066026,2.83373425054518)(10,2.78463168283138) 
};
\addlegendentry{iPotts-based};

\addplot [
color=green!50!black,
solid,
mark=x,
mark options={solid}
]
coordinates{
 (0.001,0.0896416776129664)(0.00162377673918872,0.130782005099702)(0.00263665089873036,0.171521171292151)(0.0042813323987194,0.201077582134317)(0.00695192796177561,0.237331542179857)(0.0112883789168469,0.285605323916195)(0.0183298071083244,0.357893360816685)(0.0297635144163132,0.469598850019643)(0.0483293023857175,0.646270335731638)(0.0784759970351461,0.849831768663311)(0.127427498570313,0.996468018040267)(0.206913808111479,1.10238250059031)(0.335981828628378,1.30866707331305)(0.545559478116852,1.70102549002036)(0.885866790410082,2.1898689904764)(1.43844988828766,2.79365857081443)(2.33572146909012,2.79991505277133)(3.79269019073225,2.7647742897183)(6.15848211066026,2.83373425054518)(10,2.78463168283138) 
};
\addlegendentry{Direct ADMM};

\end{semilogxaxis}
\end{tikzpicture}%
        \end{tabular}
    \caption{Reconstruction of  a jump-sparse signal using only every second frequency of the Fourier spectrum. Data is  corrupted by Gaussian noise ($\sigma = 0.05$).
The peak signal-to-noise-ratio of inverse Potts reconstructions are significantly higher than those of the TV reconstructions.
The two plots on the righthand side show the reconstruction results corresponding to the optimal regularization parameter with respect to the PSNR ($\gamma = 0.21$ for TV and $\gamma = 0.02$ for iPotts).
}
       \label{fig:pottsFourierL2}
\end{figure*}
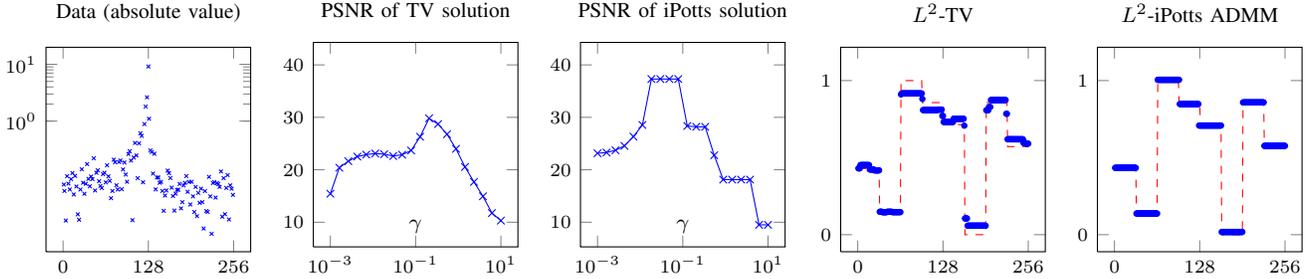

The first subproblem of the ADMM iteration \eqref{eq:ADMM} is a classical $L^2$-Potts problem of the form
\begin{equation}\label{eq:non_inverse_potts}
	P_\delta'(u) = \delta \cdot \| \nabla u \|_0 + \| u - f \|_2^2 \to \min{}
	\end{equation}
for parameter $\delta = \frac{2\gamma}{\mu}$ and data $f = v^k - \frac{\lambda^k}{\mu}.$

 For univariate data this problem can be solved fast and exactly  using  dynamic programming 
 \cite{mumford1985boundary, chambolle1995image, friedrich2008complexity, demaret2012reconstruction}.
The basic idea is that a minimizer of the Potts functional
for data $(f_1,...,f_{r})$ can be computed in polynomial time provided that minimizers of the  partial data $(f_1),$ $(f_1, f_2),$ $...,$ $ (f_1, ..., f_{r-1})$ are known.
The corresponding procedure works as follows.
We denote the respective minimizers for the partial data by $u^1,$  $u^2,$ ..., $u^{r-1}.$
In order to compute a minimizer for data $(f_1,...,f_{r}),$ 
we create a set of $r$ minimizer candidates $v^1,$ ..., $v^{r},$ each of length $r.$
These minimizer candidates are given by
\begin{equation}\label{eq:potts_candidate}
	v^\ell = (u^{\ell-1} , \underbrace{\mu_{[\ell,r]},..., \mu_{[\ell,r]}}_{\text{Length } r - \ell +1 }), \\
\end{equation}
where $u^0$ is the empty vector and 
$\mu_{[\ell, r]}$ denotes the mean value of data $f_{[\ell, r]} = (f_{\ell}, ..., f_{r}).$
 Among these candidates $v^\ell,$ one with the least Potts functional value is a minimizer for the data $f_{[1, r]}.$
{\hl 
The dynamic program for the classical Potts problem (i.e., the recursive computation of $u^n$ using \eqref{eq:potts_candidate}) can be performed in  $O(n^2)$ time and $O(n)$ space complexity \cite{friedrich2008complexity}.
 There are strategies to prune the search space  
which speed up the algorithm in practice \cite{rigaill2010pruned,
killick2012optimal}.}

For multivariate data, we cannot solve the first subproblem of our ADMM algorithm exactly in reasonable time because the classical Potts problem \eqref{eq:non_inverse_potts} is NP-hard in two dimensions \cite{boykov2001fast}. 
However, there exist well-working practical approaches based on graph cuts. We here use the max-flow/min-cut based algorithm of the library \texttt{GCOptimization 3.0} \cite{boykov2001fast, boykov2004experimental, kolmogorov2004energy}.

\subsection{Minimization of the Tikhonov subproblem}
\label{ssec:tikhonov}
The second subproblem of the ADMM iteration \eqref{eq:ADMM} is a classical Tikhonov problem
with $L^p$ data fitting of the form
\begin{equation}\label{eq:tikhonov}
	\tfrac{\mu_k}{2} \|v - w\|_2^2 + \| A v - b \|_p^p \to \mathrm{min,}
\end{equation}
where the offset vector $w$ is given by $w = u^{k+1} + \frac{\lambda^k}{\mu_k}.$ 
{\hl
The problem is convex  for 
all $p \in [1, \infty].$
Thus it can be solved efficiently using convex optimization.
We  briefly describe minimization strategies for the most relevant cases $p =1$ and $p=2.$
}

For $p = 2,$ the solution is explicitly given by the solution of the normal equation
\begin{align}\label{eq:normal_equation}
	(A^* A + \tfrac{\mu_k}{2}  \operatorname{id}) v = \tfrac{\mu_k}{2}u^{k+1} + \tfrac{1}{2}\lambda^k + A^*b.
\end{align}
Here $\operatorname{id}$ denotes the identity matrix and $A^*$ denotes the transposed of the conjugate.
As the time complexity of solving \eqref{eq:normal_equation} is $O(n^3)$ in general,
the solution of \eqref{eq:normal_equation} is the most expensive step in the ADMM iteration
since the classical univariate Potts problem is in $O(n^2).$
However, if $A^*A$ is a bandmatrix 
or if $A^*A$ can be diagonalized efficiently then the system \eqref{eq:normal_equation} can be solved fast and we are thus able to deal with large data sizes.
For instance, if
 $Ax$ describes the (circular) convolution of $x$ with some vector $h,$ i.e.,
$Ax = h * x$
then the solution of the normal equation is given by
\[
	w = \Fc^{-1}  \p{\frac{\hat{r} }{ |\hat h|^2 + \frac{\mu_k}{2}}}
\]
where $r$ denotes the right hand side of \eqref{eq:normal_equation}.

For  $p = 1$ the minimization of the Tikhonov problem \eqref{eq:tikhonov}
is more challenging because the $L^1$ data term is not differentiable.
Nevertheless, the problem can be treated by convex optimization. 
We  use the approach proposed in \cite{clason2010semismooth}.
There, the dual problem of \eqref{eq:tikhonov} is solved iteratively by a semismooth Newton method,
which converges superlinearly.
The time complexity of every iteration depends on the number of measurements
since an $m\times m$ linear system is solved in each iteration.

\section{Applications to jump-sparse recovery and numerical experiments 
}
\label{sec:experimental_jump_sparse}

In this section, we apply the inverse Potts ADMM (Algorithm \ref{alg:admm}) 
to the reconstruction of jump-sparse signals from blurred, noisy data.
We consider both reconstruction from Fourier data
and deconvolution under Gaussian, Laplacian or impulsive noise. (We refer to Appendix \ref{app:noise} for a formal description of the noise models.)
We compare the results with the minimizers of the \emph{total variation (TV)} problem  given by
\begin{equation}\label{eq:tv_min}
	\gamma \| \nabla  u \|_1 + \| A u - f\|_p^p \to \mathrm{min}.
\end{equation}
For the solution of this convex problem, 
we use the primal-dual method of \cite{chambolle2011first} with $10\,000$ iterations.

The experiments were conducted on an Apple MacBook Pro, with Intel Core 2 Duo 2.66 GHz and 8 GB RAM. Typical runtimes are between $1$ and $5$ seconds for the one-dimensional experiments, and between $5$ and $10$ minutes for two dimensions.

\subsection{Deconvolution of blurred incomplete data contaminated by Gaussian and non-Gaussian noise}\label{sec:pottsDeconv}

Here, the measurement matrix $A$ models the convolution 
with some kernel $h = (h_{-r},...,h_0, ...,h_{r})$ of non-vanishing mean.
We assume that only $m$ measurements $\{ j_1,..., j_m\},$
 $m < n,$  are given.
Hence, $A$ is a reduced $m \times n$ Toeplitz matrix of the form
\begin{equation}\label{eq:toeplitz}
	A_{j,k} = 
	\begin{cases}
	h_{k-j}, &\text{if }|k - j| \leq r \\
	0, & \text{else.}
	\end{cases}
\end{equation}
where $j = j_1,..., j_m,$ and $k = 1,...,n.$
In our experiments, $h$ is a Gaussian convolution kernel of standard deviation $6.$

In \autoref{fig:pottsDeconvL1}, data $b = A\overline{x}$ is corrupted by  Gaussian, Laplacian and impulsive noise (from top to bottom) and $m = \frac{n}{2}$ random measurements are available.
The noise variance is $\sigma = 0.05$ for Gaussian and Laplacian noise;
in the impulsive noise case, $30\%$ of the convolved signal is set to a random value between $0$ and $1$
(uniformly distributed).
For data contaminated by Gaussian noise we use the $L^2$ data term,
and for the other cases the $L^1$ data term.
In the experiment (\autoref{fig:pottsDeconvL1}) we observe that the inverse Potts ADMM algorithm performs as well as the total variation for impulsive noise.
For Gaussian and Laplacian noise, the minimizers of the total variation problem have additional plateaus as well as transitional points between the plateaus. 
In contrast, the iPotts-ADMM algorithm almost perfectly recovers the jump-sparse signal, and, in particular, the correct number of jumps.

\subsection{Reconstruction of jump-sparse signal from noisy and incomplete Fourier spectrum}

We measure an incomplete set of $m$ frequency components  of a jump-sparse signal $\overline{x} \in \R^n.$
Hence, our measurement matrix is a reduced $(m\times n)$ Fourier matrix of the form
\begin{equation*}
	A_{j,k} = \frac{1}{\sqrt{n}}e^{-2 \pi i j k / n}
\end{equation*}
where  $k = 1,...,n$ and $j$ belongs to a set of $m$ indices between $1$ and $n.$
Such reconstruction problems have been considered for example in \cite{candes2006robust, candes2006stable, yang2010fast}.
Here, we measure every second frequency component, i.e., $j = 2,4,...,n.$
We further assume that the complex valued Fourier data is corrupted by additive noise, i.e.,
\[
	b = A \overline{x} + \eta_\sigma + i \eta_\sigma'
\]
where $\eta_\sigma, \eta_\sigma'$ are $m$-dimensional vectors of i.i.d. Gaussian random variables of variance $\sigma.$

In \autoref{fig:pottsFourierL2}, we compare the performance of the inverse Potts algorithm (Algorithm \ref{alg:admm}) with that of  TV minimization \eqref{eq:tv_min}.
 We see that our method yields significantly higher  peak signal-to-noise-ratios (PSNR)
  than minimizers of the total variation problem.
 The PSNR is given by
  \begin{equation}\label{eq:psnr_def}
	\mathrm{PSNR}(x) = 10 \log_{10}\left(n \tfrac{\| \overline{x}\|_\infty^2}{  \| \overline{x} - x\|_2^2} \right)
\end{equation}
where $\overline{x}$ denotes the groundtruth.
We further observe that minimizers of the total variation problem
have small variations within the plateaus and underestimate the jump heights (\enquote{contrast reduction}).
The proposed inverse Potts ADMM algorithm reconstructs the original signal almost perfectly.

\subsection{Reconstruction and 
 segmentation of blurred images}
 
 \begin{figure}
    \setlength\figurewidth{0.15\textwidth}
    \centering
    \def\imgfoldername{experimentsRev/Potts/2D}
	\begin{subfigure}[t]{\figurewidth}
	\includegraphics[width = \figurewidth]{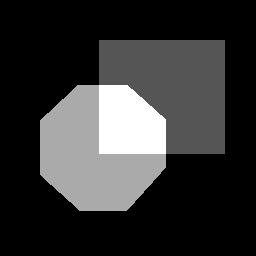} 
	\subcaption{Original}
	\end{subfigure}
	\hfill
	\begin{subfigure}[t]{\figurewidth}
	 \includegraphics[width = \figurewidth]{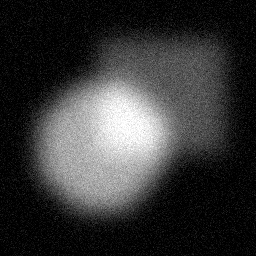}
	\subcaption{Data}
	\end{subfigure}
	\hfill
	\begin{subfigure}[t]{\figurewidth}
	 \includegraphics[width = \figurewidth]{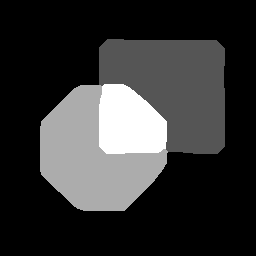}
	\subcaption{$L^2$-iPotts ADMM}
	\end{subfigure}
    \caption{
    \hl
 Deconvolution  of a geometric image  ($256 \times 256$ pixels), convolved by a Gaussian kernel of standard deviation $12$ and corrupted by  Gaussian noise $(\sigma = 0.05).$ Our method nicely removes the blur.
 The result is in particular piecewise constant as the original image.
     }
    \label{fig:pottsDeconvL2Image}
\end{figure}

{\hl
We use the inverse Potts functional in two-dimensions
for the reconstruction of cartoon-like, i.e., piecewise constant, images.
Such images serve as models in many applications, for instance in computed tomography 
\cite{ramlau2007mumford}.
In \autoref{fig:pottsDeconvL2Image}, we reconstruct a cartoon-like image from blurred and noisy data. 
Our approach recovers the piecewise constant image up to rounding off the corners. }

\begin{figure}
    \setlength\figurewidth{0.22\textwidth}
    \centering
    \def\imgfoldername{experimentsRev/Potts/2D}
    \begin{subfigure}[t]{\figurewidth}
	\includegraphics[width = \figurewidth]{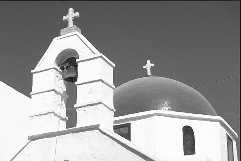} 
	\subcaption{Original}
	\end{subfigure}\hfill
	\begin{subfigure}[t]{\figurewidth}
	 \includegraphics[width = \figurewidth]{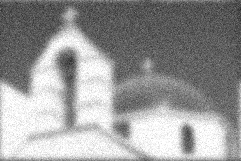}
	\subcaption{Data}
	\end{subfigure}
	\\[1ex]
	\begin{subfigure}[t]{\figurewidth}
	\includegraphics[width = \figurewidth]{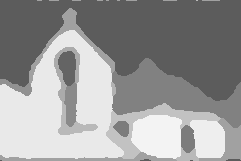}
	\subcaption{Classical $L^2$-Potts}
	\end{subfigure}\hfill
	\begin{subfigure}[t]{\figurewidth}
	 \includegraphics[width = \figurewidth]{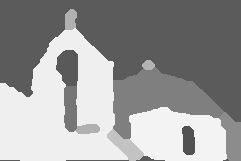}
	\subcaption{$L^2$-iPotts ADMM}
	\end{subfigure}
    \caption{Segmentation using the classical \eqref{eq:non_inverse_potts} and the inverse Potts functional
    \eqref{eq:inverse_potts} of a blurred and noisy natural image (size $241\times 161,$ image source \cite{MartinFTM01}).
    Due to the blurring, the classical Potts segmentation 
    exhibits additional segments at the boundaries, whereas segmentation with our inverse Potts ADMM detects sharp boundaries.
    }
    \label{fig:pottsDeconvL2ImageSeg}
\end{figure}

 For natural images, 
 the  Potts functional is classically used for (multi-label)  segmentation 
\cite{mumford1989optimal,boykov2001fast}. 
(The Potts problem  is sometimes called the \emph{piecewise constant Mumford-Shah problem}.) 
We see in \autoref{fig:pottsDeconvL2ImageSeg}  that the inverse Potts functional \eqref{eq:inverse_potts}, which incorporates the blurring operator $A,$ performs better than the classical Potts functional \eqref{eq:non_inverse_potts} for this task.
Here, we segment a blurred and noisy image  using  the inverse and the classical Potts functional. 
Due to the blurring, 
the segmentation using the classical Potts model introduces extra segments at the boundaries.
Minimizing the inverse Potts problem, in contrast, detects sharp boundaries without producing additional boundary segments.

\section{Applications to sparse recovery and numerical experiments
}
\label{sec:experimental_sparse}

\autoref{thm:Pottsforsparse} asserts that solutions of the inverse Potts problem associated with $A \nabla$ yield solutions of the sparsity problem 
\[
	S_\gamma(x) = \gamma \| x \|_0 + \| Ax - b\|_p^p \to \mathrm{min.}
\]
Thus, we may apply the inverse Potts ADMM (Algorithm \ref{alg:admm}) to the sparsity problem.
The corresponding method is depicted in Algorithm \ref{alg:sparsity}.
\begin{algorithm}
\small	
	  \KwIn{Data $b \in \R^m,$ model parameter $\gamma > 0,$ measurement matrix $A \in \R^{m \times n}$}
	 \KwOut{Computed result $x \in \R^n$ of the sparsity problem \eqref{eq:sparseLag} }
	 \Begin{
	 	 \SetKwHangingKw{Let}{$\leftarrow$}
	 $y $ \Let{ Solution of iPotts-ADMM (\autoref{alg:admm}) with matrix $A \nabla,$ data $f,$ and model parameter $\gamma$}
	 $x \leftarrow \nabla y$\;
	 }
	\caption{iPotts-ADMM for the sparsity problem}
	\label{alg:sparsity}
 \end{algorithm}

We compare our method (Algorithm \ref{alg:sparsity}) 
with the following  approaches  to  sparse recovery problems,
which include the state-of-the-art methods.
\begin{itemize}
\item \emph{Basis pursuit denoising (BPDN)} is the convex optimization problem
\begin{equation*}
	\gamma \, \| x\|_1 + \| Ax - b\|_p^p \to \mathrm{min}.
\end{equation*}
For the experiments, we use the toolbox \texttt{YALL1} \cite{yang2011alternating}.
\item
\emph{Iteratively reweighted $\ell^1$ minimization (IRL1)} \cite{candes2008enhancing}
  solves a sequence of constrained optimization problems
\begin{equation}\label{eq:it_rew}
	 \| x \|_{1,w} \to \mathrm{min}, \quad \text{s.t. } \| A x - b\|_2 \leq \delta,
\end{equation}
where $\| x \|_{1,w} = \sum_i w_i |x_i|$ is a weighted $\ell^1$ norm. 
The weights are initialized by $w_i = 1$ and are updated depending on the solution of the previous iteration by $w_i = \frac{1}{\epsilon + x_i}.$
We perform five iterations and choose $\epsilon = 10^{-3}.$
We use the toolbox \texttt{YALL1} \cite{yang2011alternating} for the minimization of \eqref{eq:it_rew}.
\item
 \emph{Orthogonal matching pursuit (OMP)} \cite{tropp2004greed}  greedily searches for minimizers of the constrained formulation of the $L^2$ sparsity problem
\begin{equation*}
	\min \|Ax - b\|_2^2, \quad \text{ s.t. } \| x\|_0 \leq k.
\end{equation*}
We use the implementation \texttt{OMP.m} of Stephen Becker available at Matlab's file exchange.
\item
 \emph{Iterative hard thresholding} \cite{blumensath2009iterative}
 uses surrogate functionals (forward backward splitting) for the sparsity problem.
 We here use the two variants \texttt{hard\_l0\_reg.m} \emph{(IHT-R)} and \texttt{hard\_l0\_Mterm.m} \emph{(IHT-M)} of the toolbox \texttt{sparsify 0.5}.
\item An ADMM method based on a \enquote{direct} splitting of \eqref{eq:sparseLag} which we  explain in \autoref{ssec:direct_split}.
\end{itemize}

\subsection{Reconstruction of noisy and blurred sparse signals}\label{ssec:sparse_deconv}

Our goal is to reconstruct sparse signals from noisy, blurred and incomplete measurements. 
We model this reconstruction task by \eqref{eq:sparseLag}
where $A$ is a reduced Toeplitz matrix.
In our experiments, data is blurred by a Gaussian  kernel
and $m = \frac{n}{2}$ measurements are taken.
Thus, we are in the setup of \autoref{sec:pottsDeconv}
except that now the underlying signal is sparse instead of jump-sparse.

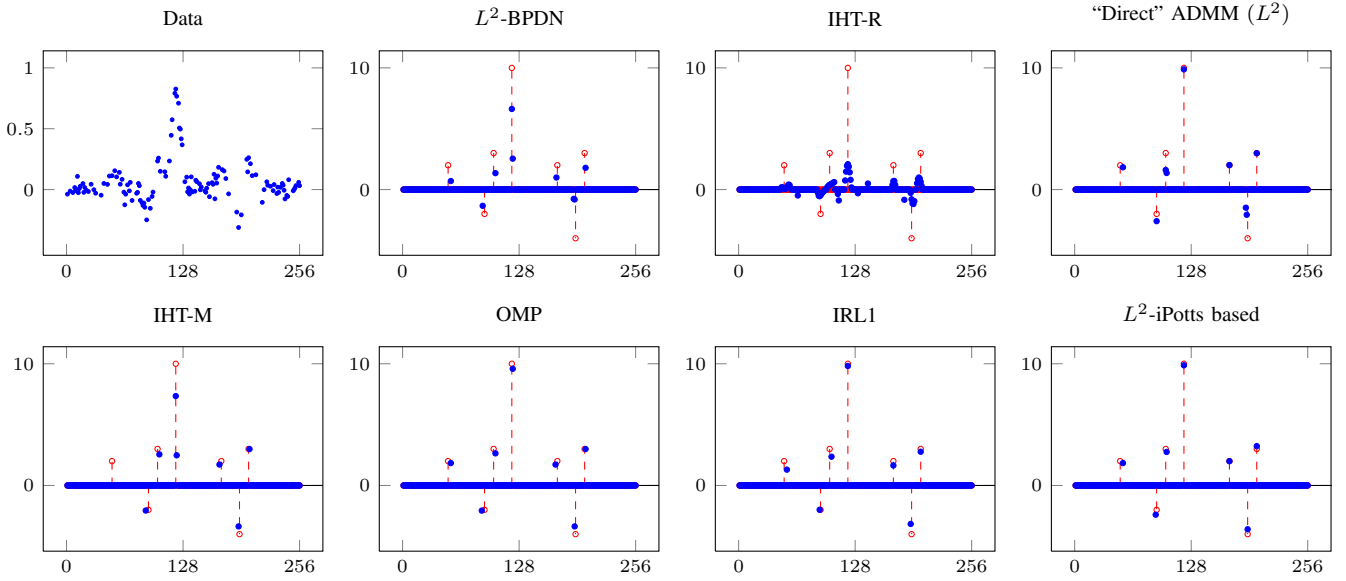
\begin{figure*}
   \setlength\figurewidth{0.205\textwidth}
    \setlength\figureheight{0.15\textwidth}
    \def\extraOptions{}
            \setlength{\tabcolsep}{-0.5ex}        
            \def\folder{experiments/L0/allPSNR/runs5/tikz/experimentGauss4/}
            \begin{tabular}{rrrr}
    \def\plottitle{Data}
    \def\nextplotname{DataGauss}
    \dataStyle
    \tikzsetnextfilename{\nextplotname}
%
%
%
\begin{tikzpicture}

\begin{semilogxaxis}[%
view={0}{90},
width=\figurewidth,
height=\figureheight,
scale only axis,
xmin=0.001, xmax=10,
xminorticks=true,
ymin=0, ymax=3,
legend style={at={(0.03,0.97)},anchor=north west,align=left},
MyAxisStyle,
title={\plottitle}]
\addplot [
color=blue,
solid,
mark=o,
mark options={solid}
]
coordinates{
 (0.001,0.100222442857688)(0.00162377673918872,0.130044761898829)(0.00263665089873036,0.163645874121925)(0.0042813323987194,0.191420362715095)(0.00695192796177561,0.224641667058485)(0.0112883789168469,0.263303122036062)(0.0183298071083244,0.317133682911229)(0.0297635144163132,0.397973172790345)(0.0483293023857175,0.527239913146773)(0.0784759970351461,0.711154350279238)(0.127427498570313,0.924204925042849)(0.206913808111479,1.04548772793932)(0.335981828628378,1.18743328122089)(0.545559478116852,1.40628893445207)(0.885866790410082,1.75743527750557)(1.43844988828766,2.45517631726791)(2.33572146909012,2.79991505277133)(3.79269019073225,2.7647742897183)(6.15848211066026,2.83373425054518)(10,2.78463168283138) 
};
\addlegendentry{iPotts-based};

\addplot [
color=green!50!black,
solid,
mark=x,
mark options={solid}
]
coordinates{
 (0.001,0.0896416776129664)(0.00162377673918872,0.130782005099702)(0.00263665089873036,0.171521171292151)(0.0042813323987194,0.201077582134317)(0.00695192796177561,0.237331542179857)(0.0112883789168469,0.285605323916195)(0.0183298071083244,0.357893360816685)(0.0297635144163132,0.469598850019643)(0.0483293023857175,0.646270335731638)(0.0784759970351461,0.849831768663311)(0.127427498570313,0.996468018040267)(0.206913808111479,1.10238250059031)(0.335981828628378,1.30866707331305)(0.545559478116852,1.70102549002036)(0.885866790410082,2.1898689904764)(1.43844988828766,2.79365857081443)(2.33572146909012,2.79991505277133)(3.79269019073225,2.7647742897183)(6.15848211066026,2.83373425054518)(10,2.78463168283138) 
};
\addlegendentry{Direct ADMM};

\end{semilogxaxis}
\end{tikzpicture}%
    &
   \sparsStyle
      \def\plottitle{$L^2$-BPDN}
    \def\nextplotname{L2iBP_YALL}
    \tikzsetnextfilename{\nextplotname}
%
%
%
\begin{tikzpicture}

\begin{semilogxaxis}[%
view={0}{90},
width=\figurewidth,
height=\figureheight,
scale only axis,
xmin=0.001, xmax=10,
xminorticks=true,
ymin=0, ymax=3,
legend style={at={(0.03,0.97)},anchor=north west,align=left},
MyAxisStyle,
title={\plottitle}]
\addplot [
color=blue,
solid,
mark=o,
mark options={solid}
]
coordinates{
 (0.001,0.100222442857688)(0.00162377673918872,0.130044761898829)(0.00263665089873036,0.163645874121925)(0.0042813323987194,0.191420362715095)(0.00695192796177561,0.224641667058485)(0.0112883789168469,0.263303122036062)(0.0183298071083244,0.317133682911229)(0.0297635144163132,0.397973172790345)(0.0483293023857175,0.527239913146773)(0.0784759970351461,0.711154350279238)(0.127427498570313,0.924204925042849)(0.206913808111479,1.04548772793932)(0.335981828628378,1.18743328122089)(0.545559478116852,1.40628893445207)(0.885866790410082,1.75743527750557)(1.43844988828766,2.45517631726791)(2.33572146909012,2.79991505277133)(3.79269019073225,2.7647742897183)(6.15848211066026,2.83373425054518)(10,2.78463168283138) 
};
\addlegendentry{iPotts-based};

\addplot [
color=green!50!black,
solid,
mark=x,
mark options={solid}
]
coordinates{
 (0.001,0.0896416776129664)(0.00162377673918872,0.130782005099702)(0.00263665089873036,0.171521171292151)(0.0042813323987194,0.201077582134317)(0.00695192796177561,0.237331542179857)(0.0112883789168469,0.285605323916195)(0.0183298071083244,0.357893360816685)(0.0297635144163132,0.469598850019643)(0.0483293023857175,0.646270335731638)(0.0784759970351461,0.849831768663311)(0.127427498570313,0.996468018040267)(0.206913808111479,1.10238250059031)(0.335981828628378,1.30866707331305)(0.545559478116852,1.70102549002036)(0.885866790410082,2.1898689904764)(1.43844988828766,2.79365857081443)(2.33572146909012,2.79991505277133)(3.79269019073225,2.7647742897183)(6.15848211066026,2.83373425054518)(10,2.78463168283138) 
};
\addlegendentry{Direct ADMM};

\end{semilogxaxis}
\end{tikzpicture}%
    &
   \sparsStyle
      \def\plottitle{IHT-R}
    \def\nextplotname{IterativeHardLag}
    \tikzsetnextfilename{\nextplotname}
%
%
%
\begin{tikzpicture}

\begin{semilogxaxis}[%
view={0}{90},
width=\figurewidth,
height=\figureheight,
scale only axis,
xmin=0.001, xmax=10,
xminorticks=true,
ymin=0, ymax=3,
legend style={at={(0.03,0.97)},anchor=north west,align=left},
MyAxisStyle,
title={\plottitle}]
\addplot [
color=blue,
solid,
mark=o,
mark options={solid}
]
coordinates{
 (0.001,0.100222442857688)(0.00162377673918872,0.130044761898829)(0.00263665089873036,0.163645874121925)(0.0042813323987194,0.191420362715095)(0.00695192796177561,0.224641667058485)(0.0112883789168469,0.263303122036062)(0.0183298071083244,0.317133682911229)(0.0297635144163132,0.397973172790345)(0.0483293023857175,0.527239913146773)(0.0784759970351461,0.711154350279238)(0.127427498570313,0.924204925042849)(0.206913808111479,1.04548772793932)(0.335981828628378,1.18743328122089)(0.545559478116852,1.40628893445207)(0.885866790410082,1.75743527750557)(1.43844988828766,2.45517631726791)(2.33572146909012,2.79991505277133)(3.79269019073225,2.7647742897183)(6.15848211066026,2.83373425054518)(10,2.78463168283138) 
};
\addlegendentry{iPotts-based};

\addplot [
color=green!50!black,
solid,
mark=x,
mark options={solid}
]
coordinates{
 (0.001,0.0896416776129664)(0.00162377673918872,0.130782005099702)(0.00263665089873036,0.171521171292151)(0.0042813323987194,0.201077582134317)(0.00695192796177561,0.237331542179857)(0.0112883789168469,0.285605323916195)(0.0183298071083244,0.357893360816685)(0.0297635144163132,0.469598850019643)(0.0483293023857175,0.646270335731638)(0.0784759970351461,0.849831768663311)(0.127427498570313,0.996468018040267)(0.206913808111479,1.10238250059031)(0.335981828628378,1.30866707331305)(0.545559478116852,1.70102549002036)(0.885866790410082,2.1898689904764)(1.43844988828766,2.79365857081443)(2.33572146909012,2.79991505277133)(3.79269019073225,2.7647742897183)(6.15848211066026,2.83373425054518)(10,2.78463168283138) 
};
\addlegendentry{Direct ADMM};

\end{semilogxaxis}
\end{tikzpicture}%
    &
   \sparsStyle
      \def\plottitle{\enquote{Direct} ADMM $(L^2)$}
    \def\nextplotname{L2iSpars_DirectADMM}
    \tikzsetnextfilename{\nextplotname}
%
%
%
\begin{tikzpicture}

\begin{semilogxaxis}[%
view={0}{90},
width=\figurewidth,
height=\figureheight,
scale only axis,
xmin=0.001, xmax=10,
xminorticks=true,
ymin=0, ymax=3,
legend style={at={(0.03,0.97)},anchor=north west,align=left},
MyAxisStyle,
title={\plottitle}]
\addplot [
color=blue,
solid,
mark=o,
mark options={solid}
]
coordinates{
 (0.001,0.100222442857688)(0.00162377673918872,0.130044761898829)(0.00263665089873036,0.163645874121925)(0.0042813323987194,0.191420362715095)(0.00695192796177561,0.224641667058485)(0.0112883789168469,0.263303122036062)(0.0183298071083244,0.317133682911229)(0.0297635144163132,0.397973172790345)(0.0483293023857175,0.527239913146773)(0.0784759970351461,0.711154350279238)(0.127427498570313,0.924204925042849)(0.206913808111479,1.04548772793932)(0.335981828628378,1.18743328122089)(0.545559478116852,1.40628893445207)(0.885866790410082,1.75743527750557)(1.43844988828766,2.45517631726791)(2.33572146909012,2.79991505277133)(3.79269019073225,2.7647742897183)(6.15848211066026,2.83373425054518)(10,2.78463168283138) 
};
\addlegendentry{iPotts-based};

\addplot [
color=green!50!black,
solid,
mark=x,
mark options={solid}
]
coordinates{
 (0.001,0.0896416776129664)(0.00162377673918872,0.130782005099702)(0.00263665089873036,0.171521171292151)(0.0042813323987194,0.201077582134317)(0.00695192796177561,0.237331542179857)(0.0112883789168469,0.285605323916195)(0.0183298071083244,0.357893360816685)(0.0297635144163132,0.469598850019643)(0.0483293023857175,0.646270335731638)(0.0784759970351461,0.849831768663311)(0.127427498570313,0.996468018040267)(0.206913808111479,1.10238250059031)(0.335981828628378,1.30866707331305)(0.545559478116852,1.70102549002036)(0.885866790410082,2.1898689904764)(1.43844988828766,2.79365857081443)(2.33572146909012,2.79991505277133)(3.79269019073225,2.7647742897183)(6.15848211066026,2.83373425054518)(10,2.78463168283138) 
};
\addlegendentry{Direct ADMM};

\end{semilogxaxis}
\end{tikzpicture}%
    \\
   \sparsStyle
      \def\plottitle{IHT-M}
    \def\nextplotname{IterativeHardCon}
    \tikzsetnextfilename{\nextplotname}
%
%
%
\begin{tikzpicture}

\begin{semilogxaxis}[%
view={0}{90},
width=\figurewidth,
height=\figureheight,
scale only axis,
xmin=0.001, xmax=10,
xminorticks=true,
ymin=0, ymax=3,
legend style={at={(0.03,0.97)},anchor=north west,align=left},
MyAxisStyle,
title={\plottitle}]
\addplot [
color=blue,
solid,
mark=o,
mark options={solid}
]
coordinates{
 (0.001,0.100222442857688)(0.00162377673918872,0.130044761898829)(0.00263665089873036,0.163645874121925)(0.0042813323987194,0.191420362715095)(0.00695192796177561,0.224641667058485)(0.0112883789168469,0.263303122036062)(0.0183298071083244,0.317133682911229)(0.0297635144163132,0.397973172790345)(0.0483293023857175,0.527239913146773)(0.0784759970351461,0.711154350279238)(0.127427498570313,0.924204925042849)(0.206913808111479,1.04548772793932)(0.335981828628378,1.18743328122089)(0.545559478116852,1.40628893445207)(0.885866790410082,1.75743527750557)(1.43844988828766,2.45517631726791)(2.33572146909012,2.79991505277133)(3.79269019073225,2.7647742897183)(6.15848211066026,2.83373425054518)(10,2.78463168283138) 
};
\addlegendentry{iPotts-based};

\addplot [
color=green!50!black,
solid,
mark=x,
mark options={solid}
]
coordinates{
 (0.001,0.0896416776129664)(0.00162377673918872,0.130782005099702)(0.00263665089873036,0.171521171292151)(0.0042813323987194,0.201077582134317)(0.00695192796177561,0.237331542179857)(0.0112883789168469,0.285605323916195)(0.0183298071083244,0.357893360816685)(0.0297635144163132,0.469598850019643)(0.0483293023857175,0.646270335731638)(0.0784759970351461,0.849831768663311)(0.127427498570313,0.996468018040267)(0.206913808111479,1.10238250059031)(0.335981828628378,1.30866707331305)(0.545559478116852,1.70102549002036)(0.885866790410082,2.1898689904764)(1.43844988828766,2.79365857081443)(2.33572146909012,2.79991505277133)(3.79269019073225,2.7647742897183)(6.15848211066026,2.83373425054518)(10,2.78463168283138) 
};
\addlegendentry{Direct ADMM};

\end{semilogxaxis}
\end{tikzpicture}%
    &
   \sparsStyle
      \def\plottitle{OMP}
    \def\nextplotname{OMP}
    \tikzsetnextfilename{\nextplotname}
%
%
%
\begin{tikzpicture}

\begin{semilogxaxis}[%
view={0}{90},
width=\figurewidth,
height=\figureheight,
scale only axis,
xmin=0.001, xmax=10,
xminorticks=true,
ymin=0, ymax=3,
legend style={at={(0.03,0.97)},anchor=north west,align=left},
MyAxisStyle,
title={\plottitle}]
\addplot [
color=blue,
solid,
mark=o,
mark options={solid}
]
coordinates{
 (0.001,0.100222442857688)(0.00162377673918872,0.130044761898829)(0.00263665089873036,0.163645874121925)(0.0042813323987194,0.191420362715095)(0.00695192796177561,0.224641667058485)(0.0112883789168469,0.263303122036062)(0.0183298071083244,0.317133682911229)(0.0297635144163132,0.397973172790345)(0.0483293023857175,0.527239913146773)(0.0784759970351461,0.711154350279238)(0.127427498570313,0.924204925042849)(0.206913808111479,1.04548772793932)(0.335981828628378,1.18743328122089)(0.545559478116852,1.40628893445207)(0.885866790410082,1.75743527750557)(1.43844988828766,2.45517631726791)(2.33572146909012,2.79991505277133)(3.79269019073225,2.7647742897183)(6.15848211066026,2.83373425054518)(10,2.78463168283138) 
};
\addlegendentry{iPotts-based};

\addplot [
color=green!50!black,
solid,
mark=x,
mark options={solid}
]
coordinates{
 (0.001,0.0896416776129664)(0.00162377673918872,0.130782005099702)(0.00263665089873036,0.171521171292151)(0.0042813323987194,0.201077582134317)(0.00695192796177561,0.237331542179857)(0.0112883789168469,0.285605323916195)(0.0183298071083244,0.357893360816685)(0.0297635144163132,0.469598850019643)(0.0483293023857175,0.646270335731638)(0.0784759970351461,0.849831768663311)(0.127427498570313,0.996468018040267)(0.206913808111479,1.10238250059031)(0.335981828628378,1.30866707331305)(0.545559478116852,1.70102549002036)(0.885866790410082,2.1898689904764)(1.43844988828766,2.79365857081443)(2.33572146909012,2.79991505277133)(3.79269019073225,2.7647742897183)(6.15848211066026,2.83373425054518)(10,2.78463168283138) 
};
\addlegendentry{Direct ADMM};

\end{semilogxaxis}
\end{tikzpicture}%
    &
   \sparsStyle
      \def\plottitle{IRL1}
    \def\nextplotname{IterativeReweight}
    \tikzsetnextfilename{\nextplotname}
%
%
%
\begin{tikzpicture}

\begin{semilogxaxis}[%
view={0}{90},
width=\figurewidth,
height=\figureheight,
scale only axis,
xmin=0.001, xmax=10,
xminorticks=true,
ymin=0, ymax=3,
legend style={at={(0.03,0.97)},anchor=north west,align=left},
MyAxisStyle,
title={\plottitle}]
\addplot [
color=blue,
solid,
mark=o,
mark options={solid}
]
coordinates{
 (0.001,0.100222442857688)(0.00162377673918872,0.130044761898829)(0.00263665089873036,0.163645874121925)(0.0042813323987194,0.191420362715095)(0.00695192796177561,0.224641667058485)(0.0112883789168469,0.263303122036062)(0.0183298071083244,0.317133682911229)(0.0297635144163132,0.397973172790345)(0.0483293023857175,0.527239913146773)(0.0784759970351461,0.711154350279238)(0.127427498570313,0.924204925042849)(0.206913808111479,1.04548772793932)(0.335981828628378,1.18743328122089)(0.545559478116852,1.40628893445207)(0.885866790410082,1.75743527750557)(1.43844988828766,2.45517631726791)(2.33572146909012,2.79991505277133)(3.79269019073225,2.7647742897183)(6.15848211066026,2.83373425054518)(10,2.78463168283138) 
};
\addlegendentry{iPotts-based};

\addplot [
color=green!50!black,
solid,
mark=x,
mark options={solid}
]
coordinates{
 (0.001,0.0896416776129664)(0.00162377673918872,0.130782005099702)(0.00263665089873036,0.171521171292151)(0.0042813323987194,0.201077582134317)(0.00695192796177561,0.237331542179857)(0.0112883789168469,0.285605323916195)(0.0183298071083244,0.357893360816685)(0.0297635144163132,0.469598850019643)(0.0483293023857175,0.646270335731638)(0.0784759970351461,0.849831768663311)(0.127427498570313,0.996468018040267)(0.206913808111479,1.10238250059031)(0.335981828628378,1.30866707331305)(0.545559478116852,1.70102549002036)(0.885866790410082,2.1898689904764)(1.43844988828766,2.79365857081443)(2.33572146909012,2.79991505277133)(3.79269019073225,2.7647742897183)(6.15848211066026,2.83373425054518)(10,2.78463168283138) 
};
\addlegendentry{Direct ADMM};

\end{semilogxaxis}
\end{tikzpicture}%
    &
   \sparsStyle
      \def\plottitle{$L^2$-iPotts based}
    \def\nextplotname{L2iSpars_PottsADMM}
    \tikzsetnextfilename{\nextplotname}
%
%
%
\begin{tikzpicture}

\begin{semilogxaxis}[%
view={0}{90},
width=\figurewidth,
height=\figureheight,
scale only axis,
xmin=0.001, xmax=10,
xminorticks=true,
ymin=0, ymax=3,
legend style={at={(0.03,0.97)},anchor=north west,align=left},
MyAxisStyle,
title={\plottitle}]
\addplot [
color=blue,
solid,
mark=o,
mark options={solid}
]
coordinates{
 (0.001,0.100222442857688)(0.00162377673918872,0.130044761898829)(0.00263665089873036,0.163645874121925)(0.0042813323987194,0.191420362715095)(0.00695192796177561,0.224641667058485)(0.0112883789168469,0.263303122036062)(0.0183298071083244,0.317133682911229)(0.0297635144163132,0.397973172790345)(0.0483293023857175,0.527239913146773)(0.0784759970351461,0.711154350279238)(0.127427498570313,0.924204925042849)(0.206913808111479,1.04548772793932)(0.335981828628378,1.18743328122089)(0.545559478116852,1.40628893445207)(0.885866790410082,1.75743527750557)(1.43844988828766,2.45517631726791)(2.33572146909012,2.79991505277133)(3.79269019073225,2.7647742897183)(6.15848211066026,2.83373425054518)(10,2.78463168283138) 
};
\addlegendentry{iPotts-based};

\addplot [
color=green!50!black,
solid,
mark=x,
mark options={solid}
]
coordinates{
 (0.001,0.0896416776129664)(0.00162377673918872,0.130782005099702)(0.00263665089873036,0.171521171292151)(0.0042813323987194,0.201077582134317)(0.00695192796177561,0.237331542179857)(0.0112883789168469,0.285605323916195)(0.0183298071083244,0.357893360816685)(0.0297635144163132,0.469598850019643)(0.0483293023857175,0.646270335731638)(0.0784759970351461,0.849831768663311)(0.127427498570313,0.996468018040267)(0.206913808111479,1.10238250059031)(0.335981828628378,1.30866707331305)(0.545559478116852,1.70102549002036)(0.885866790410082,2.1898689904764)(1.43844988828766,2.79365857081443)(2.33572146909012,2.79991505277133)(3.79269019073225,2.7647742897183)(6.15848211066026,2.83373425054518)(10,2.78463168283138) 
};
\addlegendentry{Direct ADMM};

\end{semilogxaxis}
\end{tikzpicture}%
    \end{tabular}
    \caption{ 
The original signal (dashed stem plot) is blurred by a Gaussian kernel of standard deviation $5$ and corrupted by Gaussian noise of standard deviation $\sigma = 0.05.$ We took $m = \frac{n}{2} =128$ measurements. 
Orthogonal matching pursuit (OMP), iteratively reweighted $\ell^1$ minimization (IRL1) and our iPotts based approach have the best reconstruction quality (with respect to visual inspection), followed by iterative hard thresholding (IHT-M) and the \enquote{direct} splitting 
(cf. \autoref{ssec:direct_split}).
Basis pursuit denoising (BPDN)  underestimates the heights of the spikes and the Lagrangian variant of iterative hard thresholding (IHT-R) reconstructs too many non-zero entries.
}
    \label{fig:sparsDeconvL2}
\end{figure*}

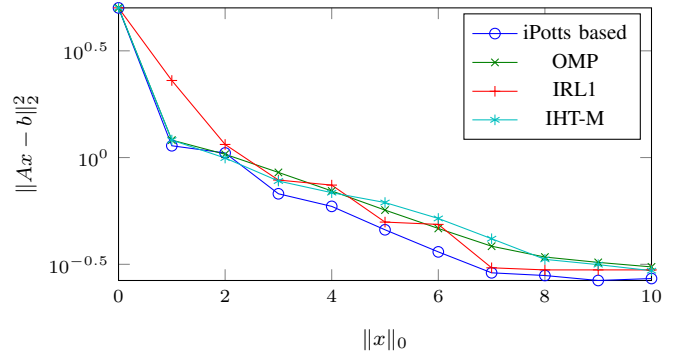
\begin{figure}
\centering
 \setlength\figurewidth{0.8\columnwidth}
 \setlength\figureheight{0.2\textwidth}
\logStyle
\def\extraOptions{}
\def\plottitle{}
%
%
%

\definecolor{mycolor1}{rgb}{0,0.75,0.75}

\begin{tikzpicture}

\begin{semilogyaxis}[%
view={0}{90},
width=\figurewidth,
height=\figureheight,
scale only axis,
xmin=0, xmax=10,
ymin=0.265552863788434, ymax=5.02648507647743,
yminorticks=true,
legend style={align=left},
MyAxisStyle,
title={\plottitle},
\extraOptions]
\addplot [
color=blue,
solid,
mark=o,
mark options={solid}
]
coordinates{
 (0,5.02648507647743)(1,1.13583425603249)(2,1.05081277940078)(3,0.676711563950321)(4,0.590688820672194)(5,0.458911278011727)(6,0.361599813402966)(7,0.288500006184523)(8,0.280195836538383)(9,0.265552863788434)(10,0.27122681860247) 
};
\addlegendentry{iPotts based};

\addplot [
color=green!50!black,
solid,
mark=x,
mark options={solid}
]
coordinates{
 (0,5.02648507647743)(1,1.20751175423534)(2,1.03623672971909)(3,0.852035112449345)(4,0.697755942786074)(5,0.567565002236034)(6,0.465549362684694)(7,0.384393749844546)(8,0.341965712388494)(9,0.3228047996263)(10,0.307510774271199) 
};
\addlegendentry{OMP};

\addplot [
color=red,
solid,
mark=+,
mark options={solid}
]
coordinates{
 (0,5.02648507647743)(1,2.29672876039551)(2,1.15194971033703)(3,0.784807085773861)(4,0.743216124701546)(5,0.498658031449839)(6,0.486233098755975)(7,0.304536442508184)(8,0.297670238626809)(9,0.297695916029016)(10,0.297694131659949) 
};
\addlegendentry{IRL1};

\addplot [
color=mycolor1,
solid,
mark=asterisk,
mark options={solid}
]
coordinates{
 (0,5.02648507647743)(1,1.20751175423534)(2,0.995115815111827)(3,0.77666772820638)(4,0.68485149093355)(5,0.615716989034622)(6,0.518482537455823)(7,0.416714472092372)(8,0.333572091192069)(9,0.315072885483047)(10,0.295202821618507) 
};
\addlegendentry{IHT-M};

\end{semilogyaxis}
\end{tikzpicture}%
\caption{
	Approximation error in dependence of the number of non-zero entries of solutions computed by state-of-the-art algorithms. Our iPotts-based method (Algorithm \ref{alg:sparsity}) yields lower approximation errors for any number of jumps.
	Here, the average values of 100 experiments are depicted.
}
\label{fig:jumps_vs_approx}
\end{figure}

Our first example is the reconstruction of blurred and incomplete data under Gaussian noise (\autoref{fig:sparsDeconvL2}).
The noise distribution suggests to employ the $L^2$ data penalty.
In the experiment, 
basis pursuit denoising (BPDN)  underestimates the height of the spikes, the Lagrangian variant of iterative hard thresholding (IHT-R) reconstructs too many non-zero entries and the \enquote{direct} splitting (\autoref{ssec:direct_split}) has to many additional non-zero entries.
Orthogonal matching pursuit (OMP), iteratively reweighted $\ell^1$ minimization (IRL1), hard thresholding (IHT-M) and the proposed iPotts-ADMM based approach 
approximate the original signal quite well; in particular, they reconstruct the precise number of non-zero entries. 
Towards a deeper comparison of these four algorithms
we quantify the reconstruction quality by looking at  
the average approximation error  $\| A x - f\|_2^2$ in dependence on the number of non-zero entries $\| x \|_0$ of a solution $x;$ cf. \autoref{fig:jumps_vs_approx}.
Here, the average values of a series of 100 runs is depicted where we used the setup of the experiment
in  \autoref{fig:sparsDeconvL2}. 
We observe that the iPotts based solutions have the least approximation errors in average.

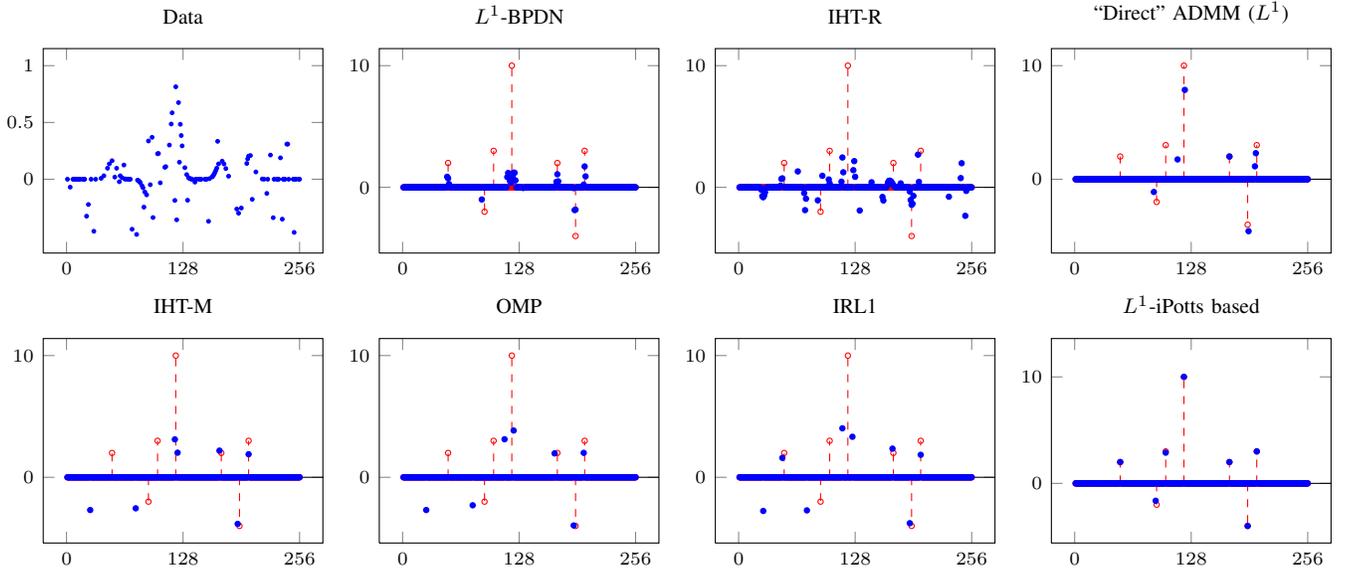
\begin{figure*}
   \setlength\figurewidth{0.205\textwidth}
    \setlength\figureheight{0.15\textwidth}
    \def\rowsep{-0.5ex}
    \def\extraOptions{}
            \setlength{\tabcolsep}{-0.5ex}        
            \def\folder{experiments/L0/allPSNR/runs5/tikz/experimentImp7/}
            \begin{tabular}{rrrr}
    \def\plottitle{Data}
    \def\nextplotname{DataImp}
    \dataStyle
    \tikzsetnextfilename{\nextplotname}
%
%
%
\begin{tikzpicture}

\begin{semilogxaxis}[%
view={0}{90},
width=\figurewidth,
height=\figureheight,
scale only axis,
xmin=0.001, xmax=10,
xminorticks=true,
ymin=0, ymax=3,
legend style={at={(0.03,0.97)},anchor=north west,align=left},
MyAxisStyle,
title={\plottitle}]
\addplot [
color=blue,
solid,
mark=o,
mark options={solid}
]
coordinates{
 (0.001,0.100222442857688)(0.00162377673918872,0.130044761898829)(0.00263665089873036,0.163645874121925)(0.0042813323987194,0.191420362715095)(0.00695192796177561,0.224641667058485)(0.0112883789168469,0.263303122036062)(0.0183298071083244,0.317133682911229)(0.0297635144163132,0.397973172790345)(0.0483293023857175,0.527239913146773)(0.0784759970351461,0.711154350279238)(0.127427498570313,0.924204925042849)(0.206913808111479,1.04548772793932)(0.335981828628378,1.18743328122089)(0.545559478116852,1.40628893445207)(0.885866790410082,1.75743527750557)(1.43844988828766,2.45517631726791)(2.33572146909012,2.79991505277133)(3.79269019073225,2.7647742897183)(6.15848211066026,2.83373425054518)(10,2.78463168283138) 
};
\addlegendentry{iPotts-based};

\addplot [
color=green!50!black,
solid,
mark=x,
mark options={solid}
]
coordinates{
 (0.001,0.0896416776129664)(0.00162377673918872,0.130782005099702)(0.00263665089873036,0.171521171292151)(0.0042813323987194,0.201077582134317)(0.00695192796177561,0.237331542179857)(0.0112883789168469,0.285605323916195)(0.0183298071083244,0.357893360816685)(0.0297635144163132,0.469598850019643)(0.0483293023857175,0.646270335731638)(0.0784759970351461,0.849831768663311)(0.127427498570313,0.996468018040267)(0.206913808111479,1.10238250059031)(0.335981828628378,1.30866707331305)(0.545559478116852,1.70102549002036)(0.885866790410082,2.1898689904764)(1.43844988828766,2.79365857081443)(2.33572146909012,2.79991505277133)(3.79269019073225,2.7647742897183)(6.15848211066026,2.83373425054518)(10,2.78463168283138) 
};
\addlegendentry{Direct ADMM};

\end{semilogxaxis}
\end{tikzpicture}%
    &
   \sparsStyle
      \def\plottitle{$L^1$-BPDN}
    \def\nextplotname{L1iBP_YALL}
    \tikzsetnextfilename{\nextplotname Imp}
%
%
%
\begin{tikzpicture}

\begin{semilogxaxis}[%
view={0}{90},
width=\figurewidth,
height=\figureheight,
scale only axis,
xmin=0.001, xmax=10,
xminorticks=true,
ymin=0, ymax=3,
legend style={at={(0.03,0.97)},anchor=north west,align=left},
MyAxisStyle,
title={\plottitle}]
\addplot [
color=blue,
solid,
mark=o,
mark options={solid}
]
coordinates{
 (0.001,0.100222442857688)(0.00162377673918872,0.130044761898829)(0.00263665089873036,0.163645874121925)(0.0042813323987194,0.191420362715095)(0.00695192796177561,0.224641667058485)(0.0112883789168469,0.263303122036062)(0.0183298071083244,0.317133682911229)(0.0297635144163132,0.397973172790345)(0.0483293023857175,0.527239913146773)(0.0784759970351461,0.711154350279238)(0.127427498570313,0.924204925042849)(0.206913808111479,1.04548772793932)(0.335981828628378,1.18743328122089)(0.545559478116852,1.40628893445207)(0.885866790410082,1.75743527750557)(1.43844988828766,2.45517631726791)(2.33572146909012,2.79991505277133)(3.79269019073225,2.7647742897183)(6.15848211066026,2.83373425054518)(10,2.78463168283138) 
};
\addlegendentry{iPotts-based};

\addplot [
color=green!50!black,
solid,
mark=x,
mark options={solid}
]
coordinates{
 (0.001,0.0896416776129664)(0.00162377673918872,0.130782005099702)(0.00263665089873036,0.171521171292151)(0.0042813323987194,0.201077582134317)(0.00695192796177561,0.237331542179857)(0.0112883789168469,0.285605323916195)(0.0183298071083244,0.357893360816685)(0.0297635144163132,0.469598850019643)(0.0483293023857175,0.646270335731638)(0.0784759970351461,0.849831768663311)(0.127427498570313,0.996468018040267)(0.206913808111479,1.10238250059031)(0.335981828628378,1.30866707331305)(0.545559478116852,1.70102549002036)(0.885866790410082,2.1898689904764)(1.43844988828766,2.79365857081443)(2.33572146909012,2.79991505277133)(3.79269019073225,2.7647742897183)(6.15848211066026,2.83373425054518)(10,2.78463168283138) 
};
\addlegendentry{Direct ADMM};

\end{semilogxaxis}
\end{tikzpicture}%
    &
   \sparsStyle
      \def\plottitle{IHT-R}
    \def\nextplotname{IterativeHardLag}
    \tikzsetnextfilename{\nextplotname Imp}
%
%
%
\begin{tikzpicture}

\begin{semilogxaxis}[%
view={0}{90},
width=\figurewidth,
height=\figureheight,
scale only axis,
xmin=0.001, xmax=10,
xminorticks=true,
ymin=0, ymax=3,
legend style={at={(0.03,0.97)},anchor=north west,align=left},
MyAxisStyle,
title={\plottitle}]
\addplot [
color=blue,
solid,
mark=o,
mark options={solid}
]
coordinates{
 (0.001,0.100222442857688)(0.00162377673918872,0.130044761898829)(0.00263665089873036,0.163645874121925)(0.0042813323987194,0.191420362715095)(0.00695192796177561,0.224641667058485)(0.0112883789168469,0.263303122036062)(0.0183298071083244,0.317133682911229)(0.0297635144163132,0.397973172790345)(0.0483293023857175,0.527239913146773)(0.0784759970351461,0.711154350279238)(0.127427498570313,0.924204925042849)(0.206913808111479,1.04548772793932)(0.335981828628378,1.18743328122089)(0.545559478116852,1.40628893445207)(0.885866790410082,1.75743527750557)(1.43844988828766,2.45517631726791)(2.33572146909012,2.79991505277133)(3.79269019073225,2.7647742897183)(6.15848211066026,2.83373425054518)(10,2.78463168283138) 
};
\addlegendentry{iPotts-based};

\addplot [
color=green!50!black,
solid,
mark=x,
mark options={solid}
]
coordinates{
 (0.001,0.0896416776129664)(0.00162377673918872,0.130782005099702)(0.00263665089873036,0.171521171292151)(0.0042813323987194,0.201077582134317)(0.00695192796177561,0.237331542179857)(0.0112883789168469,0.285605323916195)(0.0183298071083244,0.357893360816685)(0.0297635144163132,0.469598850019643)(0.0483293023857175,0.646270335731638)(0.0784759970351461,0.849831768663311)(0.127427498570313,0.996468018040267)(0.206913808111479,1.10238250059031)(0.335981828628378,1.30866707331305)(0.545559478116852,1.70102549002036)(0.885866790410082,2.1898689904764)(1.43844988828766,2.79365857081443)(2.33572146909012,2.79991505277133)(3.79269019073225,2.7647742897183)(6.15848211066026,2.83373425054518)(10,2.78463168283138) 
};
\addlegendentry{Direct ADMM};

\end{semilogxaxis}
\end{tikzpicture}%
    &
   \sparsStyle
      \def\plottitle{\enquote{Direct} ADMM ($L^1$)}
    \def\nextplotname{L1iSpars_DirectADMM}
    \tikzsetnextfilename{\nextplotname Imp}
%
%
%
\begin{tikzpicture}

\begin{semilogxaxis}[%
view={0}{90},
width=\figurewidth,
height=\figureheight,
scale only axis,
xmin=0.001, xmax=10,
xminorticks=true,
ymin=0, ymax=3,
legend style={at={(0.03,0.97)},anchor=north west,align=left},
MyAxisStyle,
title={\plottitle}]
\addplot [
color=blue,
solid,
mark=o,
mark options={solid}
]
coordinates{
 (0.001,0.100222442857688)(0.00162377673918872,0.130044761898829)(0.00263665089873036,0.163645874121925)(0.0042813323987194,0.191420362715095)(0.00695192796177561,0.224641667058485)(0.0112883789168469,0.263303122036062)(0.0183298071083244,0.317133682911229)(0.0297635144163132,0.397973172790345)(0.0483293023857175,0.527239913146773)(0.0784759970351461,0.711154350279238)(0.127427498570313,0.924204925042849)(0.206913808111479,1.04548772793932)(0.335981828628378,1.18743328122089)(0.545559478116852,1.40628893445207)(0.885866790410082,1.75743527750557)(1.43844988828766,2.45517631726791)(2.33572146909012,2.79991505277133)(3.79269019073225,2.7647742897183)(6.15848211066026,2.83373425054518)(10,2.78463168283138) 
};
\addlegendentry{iPotts-based};

\addplot [
color=green!50!black,
solid,
mark=x,
mark options={solid}
]
coordinates{
 (0.001,0.0896416776129664)(0.00162377673918872,0.130782005099702)(0.00263665089873036,0.171521171292151)(0.0042813323987194,0.201077582134317)(0.00695192796177561,0.237331542179857)(0.0112883789168469,0.285605323916195)(0.0183298071083244,0.357893360816685)(0.0297635144163132,0.469598850019643)(0.0483293023857175,0.646270335731638)(0.0784759970351461,0.849831768663311)(0.127427498570313,0.996468018040267)(0.206913808111479,1.10238250059031)(0.335981828628378,1.30866707331305)(0.545559478116852,1.70102549002036)(0.885866790410082,2.1898689904764)(1.43844988828766,2.79365857081443)(2.33572146909012,2.79991505277133)(3.79269019073225,2.7647742897183)(6.15848211066026,2.83373425054518)(10,2.78463168283138) 
};
\addlegendentry{Direct ADMM};

\end{semilogxaxis}
\end{tikzpicture}%
    \\[\rowsep]
   \sparsStyle
      \def\plottitle{IHT-M}
    \def\nextplotname{IterativeHardCon}
    \tikzsetnextfilename{\nextplotname Imp}
%
%
%
\begin{tikzpicture}

\begin{semilogxaxis}[%
view={0}{90},
width=\figurewidth,
height=\figureheight,
scale only axis,
xmin=0.001, xmax=10,
xminorticks=true,
ymin=0, ymax=3,
legend style={at={(0.03,0.97)},anchor=north west,align=left},
MyAxisStyle,
title={\plottitle}]
\addplot [
color=blue,
solid,
mark=o,
mark options={solid}
]
coordinates{
 (0.001,0.100222442857688)(0.00162377673918872,0.130044761898829)(0.00263665089873036,0.163645874121925)(0.0042813323987194,0.191420362715095)(0.00695192796177561,0.224641667058485)(0.0112883789168469,0.263303122036062)(0.0183298071083244,0.317133682911229)(0.0297635144163132,0.397973172790345)(0.0483293023857175,0.527239913146773)(0.0784759970351461,0.711154350279238)(0.127427498570313,0.924204925042849)(0.206913808111479,1.04548772793932)(0.335981828628378,1.18743328122089)(0.545559478116852,1.40628893445207)(0.885866790410082,1.75743527750557)(1.43844988828766,2.45517631726791)(2.33572146909012,2.79991505277133)(3.79269019073225,2.7647742897183)(6.15848211066026,2.83373425054518)(10,2.78463168283138) 
};
\addlegendentry{iPotts-based};

\addplot [
color=green!50!black,
solid,
mark=x,
mark options={solid}
]
coordinates{
 (0.001,0.0896416776129664)(0.00162377673918872,0.130782005099702)(0.00263665089873036,0.171521171292151)(0.0042813323987194,0.201077582134317)(0.00695192796177561,0.237331542179857)(0.0112883789168469,0.285605323916195)(0.0183298071083244,0.357893360816685)(0.0297635144163132,0.469598850019643)(0.0483293023857175,0.646270335731638)(0.0784759970351461,0.849831768663311)(0.127427498570313,0.996468018040267)(0.206913808111479,1.10238250059031)(0.335981828628378,1.30866707331305)(0.545559478116852,1.70102549002036)(0.885866790410082,2.1898689904764)(1.43844988828766,2.79365857081443)(2.33572146909012,2.79991505277133)(3.79269019073225,2.7647742897183)(6.15848211066026,2.83373425054518)(10,2.78463168283138) 
};
\addlegendentry{Direct ADMM};

\end{semilogxaxis}
\end{tikzpicture}%
    &
   \sparsStyle
      \def\plottitle{OMP}
    \def\nextplotname{OMP}
    \tikzsetnextfilename{\nextplotname Imp}
%
%
%
\begin{tikzpicture}

\begin{semilogxaxis}[%
view={0}{90},
width=\figurewidth,
height=\figureheight,
scale only axis,
xmin=0.001, xmax=10,
xminorticks=true,
ymin=0, ymax=3,
legend style={at={(0.03,0.97)},anchor=north west,align=left},
MyAxisStyle,
title={\plottitle}]
\addplot [
color=blue,
solid,
mark=o,
mark options={solid}
]
coordinates{
 (0.001,0.100222442857688)(0.00162377673918872,0.130044761898829)(0.00263665089873036,0.163645874121925)(0.0042813323987194,0.191420362715095)(0.00695192796177561,0.224641667058485)(0.0112883789168469,0.263303122036062)(0.0183298071083244,0.317133682911229)(0.0297635144163132,0.397973172790345)(0.0483293023857175,0.527239913146773)(0.0784759970351461,0.711154350279238)(0.127427498570313,0.924204925042849)(0.206913808111479,1.04548772793932)(0.335981828628378,1.18743328122089)(0.545559478116852,1.40628893445207)(0.885866790410082,1.75743527750557)(1.43844988828766,2.45517631726791)(2.33572146909012,2.79991505277133)(3.79269019073225,2.7647742897183)(6.15848211066026,2.83373425054518)(10,2.78463168283138) 
};
\addlegendentry{iPotts-based};

\addplot [
color=green!50!black,
solid,
mark=x,
mark options={solid}
]
coordinates{
 (0.001,0.0896416776129664)(0.00162377673918872,0.130782005099702)(0.00263665089873036,0.171521171292151)(0.0042813323987194,0.201077582134317)(0.00695192796177561,0.237331542179857)(0.0112883789168469,0.285605323916195)(0.0183298071083244,0.357893360816685)(0.0297635144163132,0.469598850019643)(0.0483293023857175,0.646270335731638)(0.0784759970351461,0.849831768663311)(0.127427498570313,0.996468018040267)(0.206913808111479,1.10238250059031)(0.335981828628378,1.30866707331305)(0.545559478116852,1.70102549002036)(0.885866790410082,2.1898689904764)(1.43844988828766,2.79365857081443)(2.33572146909012,2.79991505277133)(3.79269019073225,2.7647742897183)(6.15848211066026,2.83373425054518)(10,2.78463168283138) 
};
\addlegendentry{Direct ADMM};

\end{semilogxaxis}
\end{tikzpicture}%
    &
   \sparsStyle
      \def\plottitle{IRL1}
    \def\nextplotname{IterativeReweight}
    \tikzsetnextfilename{\nextplotname Imp}
%
%
%
\begin{tikzpicture}

\begin{semilogxaxis}[%
view={0}{90},
width=\figurewidth,
height=\figureheight,
scale only axis,
xmin=0.001, xmax=10,
xminorticks=true,
ymin=0, ymax=3,
legend style={at={(0.03,0.97)},anchor=north west,align=left},
MyAxisStyle,
title={\plottitle}]
\addplot [
color=blue,
solid,
mark=o,
mark options={solid}
]
coordinates{
 (0.001,0.100222442857688)(0.00162377673918872,0.130044761898829)(0.00263665089873036,0.163645874121925)(0.0042813323987194,0.191420362715095)(0.00695192796177561,0.224641667058485)(0.0112883789168469,0.263303122036062)(0.0183298071083244,0.317133682911229)(0.0297635144163132,0.397973172790345)(0.0483293023857175,0.527239913146773)(0.0784759970351461,0.711154350279238)(0.127427498570313,0.924204925042849)(0.206913808111479,1.04548772793932)(0.335981828628378,1.18743328122089)(0.545559478116852,1.40628893445207)(0.885866790410082,1.75743527750557)(1.43844988828766,2.45517631726791)(2.33572146909012,2.79991505277133)(3.79269019073225,2.7647742897183)(6.15848211066026,2.83373425054518)(10,2.78463168283138) 
};
\addlegendentry{iPotts-based};

\addplot [
color=green!50!black,
solid,
mark=x,
mark options={solid}
]
coordinates{
 (0.001,0.0896416776129664)(0.00162377673918872,0.130782005099702)(0.00263665089873036,0.171521171292151)(0.0042813323987194,0.201077582134317)(0.00695192796177561,0.237331542179857)(0.0112883789168469,0.285605323916195)(0.0183298071083244,0.357893360816685)(0.0297635144163132,0.469598850019643)(0.0483293023857175,0.646270335731638)(0.0784759970351461,0.849831768663311)(0.127427498570313,0.996468018040267)(0.206913808111479,1.10238250059031)(0.335981828628378,1.30866707331305)(0.545559478116852,1.70102549002036)(0.885866790410082,2.1898689904764)(1.43844988828766,2.79365857081443)(2.33572146909012,2.79991505277133)(3.79269019073225,2.7647742897183)(6.15848211066026,2.83373425054518)(10,2.78463168283138) 
};
\addlegendentry{Direct ADMM};

\end{semilogxaxis}
\end{tikzpicture}%
    &
   \sparsStyle
      \def\plottitle{$L^1$-iPotts based}
    \def\nextplotname{L1iSpars_PottsADMM}
    \tikzsetnextfilename{\nextplotname Imp}
%
%
%
\begin{tikzpicture}

\begin{semilogxaxis}[%
view={0}{90},
width=\figurewidth,
height=\figureheight,
scale only axis,
xmin=0.001, xmax=10,
xminorticks=true,
ymin=0, ymax=3,
legend style={at={(0.03,0.97)},anchor=north west,align=left},
MyAxisStyle,
title={\plottitle}]
\addplot [
color=blue,
solid,
mark=o,
mark options={solid}
]
coordinates{
 (0.001,0.100222442857688)(0.00162377673918872,0.130044761898829)(0.00263665089873036,0.163645874121925)(0.0042813323987194,0.191420362715095)(0.00695192796177561,0.224641667058485)(0.0112883789168469,0.263303122036062)(0.0183298071083244,0.317133682911229)(0.0297635144163132,0.397973172790345)(0.0483293023857175,0.527239913146773)(0.0784759970351461,0.711154350279238)(0.127427498570313,0.924204925042849)(0.206913808111479,1.04548772793932)(0.335981828628378,1.18743328122089)(0.545559478116852,1.40628893445207)(0.885866790410082,1.75743527750557)(1.43844988828766,2.45517631726791)(2.33572146909012,2.79991505277133)(3.79269019073225,2.7647742897183)(6.15848211066026,2.83373425054518)(10,2.78463168283138) 
};
\addlegendentry{iPotts-based};

\addplot [
color=green!50!black,
solid,
mark=x,
mark options={solid}
]
coordinates{
 (0.001,0.0896416776129664)(0.00162377673918872,0.130782005099702)(0.00263665089873036,0.171521171292151)(0.0042813323987194,0.201077582134317)(0.00695192796177561,0.237331542179857)(0.0112883789168469,0.285605323916195)(0.0183298071083244,0.357893360816685)(0.0297635144163132,0.469598850019643)(0.0483293023857175,0.646270335731638)(0.0784759970351461,0.849831768663311)(0.127427498570313,0.996468018040267)(0.206913808111479,1.10238250059031)(0.335981828628378,1.30866707331305)(0.545559478116852,1.70102549002036)(0.885866790410082,2.1898689904764)(1.43844988828766,2.79365857081443)(2.33572146909012,2.79991505277133)(3.79269019073225,2.7647742897183)(6.15848211066026,2.83373425054518)(10,2.78463168283138) 
};
\addlegendentry{Direct ADMM};

\end{semilogxaxis}
\end{tikzpicture}%
    \end{tabular}
    \caption{The same setup as in \autoref{fig:sparsDeconvL2} 
    replacing Gaussian noise by impulsive noise ($25\%$ of data were set to a random value between $-0.5$ and $0.5$).
For the direct splitting (cf. \autoref{ssec:direct_split}), basis pursuit (BPDN),
and the iPotts based method we use  $L^1$ data terms.
We see that the proposed $L^1$-iPotts based algorithm performs significantly better than the other methods.
It is able to recover the original signal almost perfectly.
}
    \label{fig:sparsDeconvL1}
\end{figure*}

In  \autoref{fig:sparsDeconvL1}, we drive the same experiment as in \autoref{fig:sparsDeconvL2}
replacing Gaussian noise by impulsive noise.
Due to this noise model, we employ the $L^1$ data term for our iPotts-based algorithm.
For the other methods we also use the $L^1$ variant whenever it is available;
to the best of our knowledge, this is the case for basis pursuit denoising and
the direct splitting (\autoref{ssec:direct_split}).
We observe that the proposed algorithm yields an almost perfect reconstruction also in presence of impulsive noise and that it performs significantly better than the other methods in this case.

\subsection{Comparison with a \enquote{direct} ADMM approach to the sparsity problem}\label{ssec:direct_split}

In analogy to \eqref{eq:potts_split}, we consider the consensus form of the sparsity problem 
 \begin{equation}\label{eq:spars_consensus}
\gamma \| u \|_0 + \| A v - b \|_p^p \to \text{min,}\quad 
\text{s.t. }  u -  v = 0.
\end{equation}
This leads to the augmented Lagrangian
 \begin{equation}\label{eq:lagrangian_spars}
\begin{split}
 	 &  \gamma \| u \|_0 + \langle \lambda , u - v \rangle  \\
	& +  \tfrac{\mu}{2}  \| u-v  \|_2^2 + \| A v - b \|_p^p \to \mathrm{min.}
	\end{split}
\end{equation}
Proceeding as in \autoref{sec:Potts_ADMM} we obtain a \enquote{direct} ADMM algorithm for the sparsity problem.
This algorithm is given by replacing $\| \nabla u \|_0$ by $\| u \|_0$ in the first line of \eqref{eq:ADMM}. 
This leads to alternately solving a hard thresholding problem (instead of a Potts problem) and a classical Tikhonov problem associated with matrix $A.$

The difference between the \enquote{direct} ADMM approach and
 our iPotts-ADMM based method (Algorithm \ref{alg:sparsity}) is that they are based on different augmented Lagrangians.
Indeed, when applying the iPotts-ADMM to the sparsity problem,
we consider the inverse Potts problem associated with $A \nabla$ instead of $A.$
Then, the augmented Lagrangian of the corresponding problem is obtained by replacing 
$A$ by $A\nabla$ in \eqref{eq:lagrangian}.
With the substitutions $\nabla u = u'$ and $\nabla v = v',$ 
equation \eqref{eq:lagrangian} reads
 \begin{equation}\label{eq:lagrangian_iPotts}
\begin{split}
 	\, &  \gamma \| u' \|_0 + \langle \lambda , \nabla^+(u' - v') \rangle  \\
	& +  \tfrac{\mu}{2}  \| \nabla^+(u'-v')  \|_2^2 + \| A v' - b \|_p^p \to \mathrm{min.}
	\end{split}
\end{equation}
Comparing \eqref{eq:lagrangian_spars} and \eqref{eq:lagrangian_iPotts},
we see that the direct method couples $u$ and $v$ directly whereas
the iPotts based method involves the  antiderivatives of $u$ and $v.$

\begin{figure}
\centering
\setlength\figurewidth{0.4\textwidth}
    \setlength\figureheight{0.15\textwidth}
      \pgfkeys{/pgfplots/MyAxisStyle/.style={
        font=\footnotesize,
        tick label style = {font=\scriptsize},
        enlargelimits, 
        y label style={at={(0.1,0.5)}},
        xlabel = {$\gamma$},
        ylabel = {Average energy $S_\gamma(x)$}
	}}
    \def\plottitle{}
    \def\nextplotname{demoADMMEnergy}
    \tikzsetnextfilename{\nextplotname}
%
%
%
\begin{tikzpicture}

\begin{semilogxaxis}[%
view={0}{90},
width=\figurewidth,
height=\figureheight,
scale only axis,
xmin=0.001, xmax=10,
xminorticks=true,
ymin=0, ymax=3,
legend style={at={(0.03,0.97)},anchor=north west,align=left},
MyAxisStyle,
title={\plottitle}]
\addplot [
color=blue,
solid,
mark=o,
mark options={solid}
]
coordinates{
 (0.001,0.100222442857688)(0.00162377673918872,0.130044761898829)(0.00263665089873036,0.163645874121925)(0.0042813323987194,0.191420362715095)(0.00695192796177561,0.224641667058485)(0.0112883789168469,0.263303122036062)(0.0183298071083244,0.317133682911229)(0.0297635144163132,0.397973172790345)(0.0483293023857175,0.527239913146773)(0.0784759970351461,0.711154350279238)(0.127427498570313,0.924204925042849)(0.206913808111479,1.04548772793932)(0.335981828628378,1.18743328122089)(0.545559478116852,1.40628893445207)(0.885866790410082,1.75743527750557)(1.43844988828766,2.45517631726791)(2.33572146909012,2.79991505277133)(3.79269019073225,2.7647742897183)(6.15848211066026,2.83373425054518)(10,2.78463168283138) 
};
\addlegendentry{iPotts-based};

\addplot [
color=green!50!black,
solid,
mark=x,
mark options={solid}
]
coordinates{
 (0.001,0.0896416776129664)(0.00162377673918872,0.130782005099702)(0.00263665089873036,0.171521171292151)(0.0042813323987194,0.201077582134317)(0.00695192796177561,0.237331542179857)(0.0112883789168469,0.285605323916195)(0.0183298071083244,0.357893360816685)(0.0297635144163132,0.469598850019643)(0.0483293023857175,0.646270335731638)(0.0784759970351461,0.849831768663311)(0.127427498570313,0.996468018040267)(0.206913808111479,1.10238250059031)(0.335981828628378,1.30866707331305)(0.545559478116852,1.70102549002036)(0.885866790410082,2.1898689904764)(1.43844988828766,2.79365857081443)(2.33572146909012,2.79991505277133)(3.79269019073225,2.7647742897183)(6.15848211066026,2.83373425054518)(10,2.78463168283138) 
};
\addlegendentry{Direct ADMM};

\end{semilogxaxis}
\end{tikzpicture}%
    \caption{
    	Final total energy  $S_\gamma(x)$  of the sparsity problem using
	the iPotts based method (Algorithm \ref{alg:sparsity}) and the \enquote{direct} ADMM method (\autoref{ssec:direct_split}) for different parameters $\gamma.$ 
	The iPotts-ADMM algorithm reaches lower energies over the whole parameter range.
	This indicates the superiority of the first method.
	Computed values are averages over 50 experiments.
    }
    \label{fig:direct_splitting}
\end{figure}
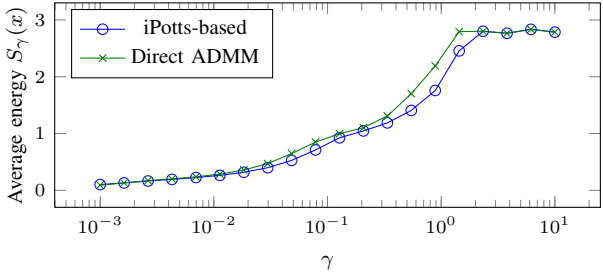

From the experiments (Figures \ref{fig:sparsDeconvL2}, \ref{fig:sparsDeconvL1}, 
\ref{fig:direct_splitting})
we conclude that the iPotts-based method (Algorithm \ref{alg:sparsity}) is advantageous over the direct ADMM. 
In particular, the solutions of 
the iPotts-based method 
have lower energy than the \enquote{direct} method for the whole range of parameters $\gamma$;
cf. \autoref{fig:direct_splitting}.

\subsection{Sparse image recovery}

\begin{figure}
    \setlength\figurewidth{0.14\textwidth}
    \footnotesize
    \centering
    \def\imgfoldername{experimentsRev/L0/2D}
    \def\rowsep{-0.5ex}
	\begin{subfigure}[t]{\figurewidth}
	\includegraphics[width = \figurewidth]{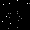} 
	\subcaption{Original}
	\end{subfigure}\hfill
	\begin{subfigure}[t]{\figurewidth}
	 \includegraphics[width = \figurewidth]{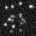}
	\subcaption{Blurred, noisy data}
	\end{subfigure}\hfill
	\begin{subfigure}[t]{\figurewidth}
	 \includegraphics[width = \figurewidth]{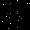}
	\subcaption{iPotts based reconstruction}
	\end{subfigure}
    \caption{Sparse image of size $30 \times 30$ 
    blurred by a $7 \times 7$ Gaussian kernel (standard deviation $1$) and corrupted by Gaussian noise ($\sigma = 0.05$).
    Data consists of  50 \% randomly selected pixels of the middle image.
The iPotts based method recovers almost all spikes correctly.
(For visualization purposes, the contrast of the middle image was increased.)}
    \label{fig:sparsDeconvL2Image}
\end{figure}

We also use our method to reconstruct  sparse images.
One may think of an image of small particles or of an astronomic image.
We apply our procedure  to images by reshaping the image to a vector and adapting the matrix $A$ accordingly. 
 \autoref{fig:sparsDeconvL2Image} shows the deconvolution of a sparse image using our iPotts-ADMM based method (Algorithm \ref{alg:sparsity}). In the experiment we see that almost all spikes are recovered while only few false positives are reconstructed.

\section{Proofs}
\label{sec:proofs}

Here we provide the proofs of the theorems stated in the course of this paper.

\subsection{Existence of Minimizers}\label{ssec:existMin}

We start out showing \autoref{thm:minimizer_existence} which asserts that the inverse Potts problem \eqref{eq:inverse_potts}
has a minimizer.

\begin{proof}[Proof of  \autoref{thm:minimizer_existence}]
In order to deal with the general case of a (possibly) singular matrix $A$ we decompose the domain into $\ker A$
and a corresponding algebraic complement $U.$ This means that $U + \ker A = \R^n$
(or $\C^n)$ and $U \cap \ker A = \{0\}$.
(In the following we proceed without drawing attention to $\C^n$ when writing $\R^n,$ but the arguments work for the complex case as well.)
For $x \in \R^n$ we frequently use the decomposition $x = u + v,$ where $u=Q_Ux$ is the projection $Q_U$ of $x$
to $U,$ and $v$ is the corresponding projection onto $\ker A.$

The matrix $A$ restricted to the subspace $U$ is invertible, and since we are in finite dimensional space,
there is a positive constant $c$ such that
$$\|Au\| \geq c \|u\| \text{ for any }u \in U.$$
(Due to the finite dimension all norms are equivalent and the above inequality holds for any norm.)
As a consequence, whenever, for a sequence $u_k$ in U, the norm $\|u_k\|$ tends to
$\infty,$ the inverse Potts functional $P_\gamma(u_k)$ defined by \eqref{eq:inverse_potts} tends to $\infty$ as well.
Therefore, for any sequence of vectors $x_k$ in $\R^n$ (not only in U,) we obtain the implication:
\begin{multline}\label{eq:EqBdProof}
    P_\gamma(x_k) \text{ is bounded } \\
    \implies  u_k = Q_U x_k \text{ has a converging subsequence.}
\end{multline}
This is a consequence of $Ax_k = AQ_Ux_k.$

Our next preparatory step introduces the mapping $s$ on $U$ which assigns to each $u \in U$ the minimal number of jumps
of all vectors in $u + \ker A,$ i.e.,
\begin{align*}
s(u) = \min_{v \in u+\ker A} \|\nabla v\|_0.
\end{align*}
We show that this mapping $s$ is lower semicontinuous which, in our context, means that the preimages of
the sets $\{0,\ldots, k \}$ are closed for all $k \in \N.$
To see this, we first observe that the set $M_k$ of all vectors in $\R^n$ with at most $k$ jumps
is structurally a finite union of vector spaces (of dimension $k+1$.) More precisely,
\[
M_k = \{x \in \R^n: \|\nabla x\|_0 \leq k \} = \bigcup_{J \subset \{1,\ldots,n-1\},|J|= k} X_J,
\]
where $X_J$ are those vectors whose jump sets are contained in $J \subset \{1,\ldots,n-1\}.$
Furthermore, a vector $u \in U$ has the property $s(u) \leq k$ if and only if there is a
vector $x \in M_k$ (i.e., with at most $k$ jumps) such that $Q_Ux=u.$
Summing up,
\begin{equation}\label{eq:EqToDeriveLSC}
 s^{-1} (\{0,\ldots, k \})  = Q_U(M_k) = \bigcup_{J \subset \{1,\ldots,n-1\},|J|= k} Q_U(X_J).
\end{equation}
We discuss the right hand side of \eqref{eq:EqToDeriveLSC} to see the lower semicontinuity of $s.$
Each $\hl Q_U(X_J)$ is a finite dimensional linear subspace and thus closed; so as a finite union
of closed sets the right hand side of \eqref{eq:EqToDeriveLSC} is closed.
Therefore the left hand side of \eqref{eq:EqToDeriveLSC} is closed which
by definition implies the lower semicontinuity of $s.$

Now we can show the assertion of the theorem. We consider a sequence $x_k$ such that the values $P_\gamma(x_k)$
of the inverse Potts functional $P_\gamma$ tend to an infimum, i.e.,
\[
  \lim_{k \to \infty}     P_\gamma(x_k) = \inf_{x \in \R^n} P_\gamma(x).
\]
For every member of the sequence, we write $x_k = u_k + v_k$ with $u_k \in U$ and $v_k \in \ker A.$
By \eqref{eq:EqBdProof} we find a subsequence $x_{k_l}$ such that $u_{k_l} = Q_U x_{k_l}$ converges to some $u \in U.$
Since $P_\gamma(x_{k_l})$ converges and $Au_{k_l} = Ax_{k_l}$  we have that
\begin{multline*}
  | \|\nabla x_{k_l}\|_0 - \|\nabla x_{k_r}\|_0 | \leq |P_\gamma(x_{k_l}) - P_\gamma(x_{k_r})| \\
  + | \|Au_{k_l}-b \|_p^p - \|Au_{k_r}-b \|_p^p|
  \to 0  \text{ as } l,r \to \infty.
\end{multline*}
This means that, for sufficiently large $l$, the number of jumps $\|\nabla x_{k_l}\|_0$ becomes constant; let us denote this constant by $j$.
As a consequence, $s(u_{k_l}) \leq \|\nabla x_{k_l}\|_0 = j,$ and thus, by the lower semicontinuity of $s$, $s(u) \leq j.$
Hence, by the definition of $s,$ there is a vector $x^\ast \in u+ \ker A$ such that the number of jumps of $x^\ast$ is smaller than or equal to $j.$
Then,
\begin{align*}
P_\gamma(x^\ast) &=  \gamma \|\nabla x^\ast\|_0 + \|Ax^\ast-b \|_p^p
   = \gamma \|\nabla x^\ast\|_0 + \|Au-b \|_p^p \\
 &\leq j + \lim_l \|Au_{k_l}-b \|_p^p = \lim_l P_\gamma(x_{k_l})
\end{align*}
which shows that $x^\ast$ is a minimizer as desired.
\end{proof}

Next we show \autoref{thm:minimizer_not_existence_continuous_time} which states that
the continuous-time analogue of \autoref{thm:minimizer_existence} is wrong. We give
counterexamples, i.e., we find bounded operators $A$ and data $f$ such that
the continuous-time inverse Potts functional
\[
    P_\gamma(u) =
    \gamma \cdot \| \nabla  u \|_0 + \| A u -f \|_p^p,  \\
\]
if $u$ is a piecewise constant function on the interval $[0,1]$, and $P_\gamma(u) = \infty$ otherwise,
has no minimizer.

\begin{proof}[Proof of \autoref{thm:minimizer_not_existence_continuous_time}]

We consider a positive function $g \in L^p[0,1]$ with total mass $1$ which is supported in the interval
$[ \tfrac{1}{2} - \varepsilon, \tfrac{1}{2} + \varepsilon]$ with positive $\varepsilon<\tfrac{1}{8}.$
We use the symbol $\tilde{g}$ for its left-shift by $\tfrac{1}{2}.$
Our counterexamples are the (cyclic) convolution operators with functions $\tilde{g}$ as above, i.e.,
operators $A$ defined by $Au = \tilde{g} * u,$ and the data given by $f=g.$

We claim that, for Potts parameter $\gamma$ with $\gamma$ $<\gamma_0$ (defined in \eqref{eq:gamma0} below,)
\begin{equation}\label{eq:MiniNotExEq}
\inf_v P_\gamma(v) =  2 \gamma < P_\gamma(u) \quad \text{ for all }u.
\end{equation}
This means that there is no minimizer in that case and thus shows the assertion of the theorem.
In order to show the equality in \eqref{eq:MiniNotExEq}, we consider
the sequence of characteristic functions  $u_n =  \frac{n}{2} 1_{[\frac{1}{2}-\frac{1}{n},\frac{1}{2} + \frac{1}{n}]}.$
We have that $\| \nabla u_n \|_0 = 2$  and $\|A u_n -  f\|  \to 0.$ Thus, $P_\gamma(u_n) \to  2 \gamma.$
This yields $\inf_v P_\gamma(v) \leq 2 \gamma$. It remains to show the inequality in \eqref{eq:MiniNotExEq} (which in turn
implies the equality in \eqref{eq:MiniNotExEq}.) To this end, we have to consider the set of functions $u$ with at most
one jump and find $\gamma>0$ such that $d(u) = \| A u -f \|_p^p > 2 \gamma$ for all such $u.$
The set $B = \{x \in [0,1]: f(x)\geq 2\}$ has positive Lebesgue measure $\lambda(B)$ since $f$ has total mass $1$
and is supported on an interval of length bounded by $\tfrac{1}{4}.$ If $u<1,$ then $d(u)\geq \lambda(B).$
So in order to obtain $d(u) < \lambda(B),$ we need that $u \geq 1$ either to the left or to the right of its (sole) jump location.
Then $Au \geq 1$ on at least one of the intervals $[\varepsilon, \tfrac{1}{2} - \varepsilon]$
and $[\tfrac{1}{2} + \varepsilon, 1 - \varepsilon ].$
Both of these intervals have length $\tfrac{1}{2}-2 \varepsilon,$ and, on both intervals,  $f=0.$
Therefore, if $d(u) < \lambda(B),$ we necessarily have $d(u) \geq \tfrac{1}{2}-2 \varepsilon > \tfrac{1}{4}.$
Then, for any $u$ with at most one jump,
\begin{equation}\label{eq:gamma0}
d(u) > \min(\lambda(B), \tfrac{1}{4}) =: 2 \gamma_0.
\end{equation}
If $u$ has two or more jumps then trivially $P_\gamma(u) > 2 \gamma.$
Together, this implies that, for any $\gamma$ with $\gamma<\gamma_0$, the inverse Potts functional $P_\gamma$ fulfills
$P_\gamma(u) > 2\gamma$ for all $u \in L^p[0,1].$ This shows \eqref{eq:MiniNotExEq} which completes the proof.
\end{proof}

\subsection{Relations to sparsity}
\label{ssec:RelSparse}

We first prove \autoref{thm:Pottsforsparse} which shows how to transform a sparsity problem
into a jump-sparsity problem.

\begin{proof}[Proof of \autoref{thm:Pottsforsparse}]
For $x^*$ satisfying \eqref{eq:inverse_potts2}, we define
 $u^* = \nabla x^*$.  Towards a contradiction we assume that there is
  $u \in \R^n$  such that
$    \gamma \|u\|_0 + \|Au-b\|_p^p <  \gamma \|u^*\|_0 + \|Au^*-b\|_p^p$.
Then, for $u$,
there is $x \in \R^{n+1}$ such that $u = \nabla x$. Then,
\begin{align*}
   \gamma \|\nabla x\|_0 + \|A \nabla x -b\|_p^p &< \gamma \|u^*\|_0 + \|A u^* -b\|_p^p \\
   &= \gamma \|\nabla x^*\|_0 + \|A \nabla x^* -b\|_p^p.
\end{align*}
which is a contradiction.
\end{proof}

For $p=2$ we show a converse statement. It is formulated as \autoref{thm:sparseforPotts} and proved next.
In its proof we make use of the decomposition of $\R^n$ into the orthogonal direct sum $\R^n= V \oplus \R e$,
where $e$ denoted the constant vector $(1,\hdots,1)^T$ and $V$ is the linear
space of vectors with zero mean. Observing that the linear operator $\nabla$ is bijective from the linear space $V$
to $\R^{n-1}$, we use the symbol $\nabla^+$ for the mapping $\R^{n-1} \to V,$
\begin{equation} \label{eq:defNablaPlus}
   \nabla^+ =  (\nabla|_V)^{-1},
\end{equation}
for the inverse of the mapping $\nabla$ restricted to the subspace $V$.

\begin{proof}[Proof of \autoref{thm:sparseforPotts}]
We consider the inverse Potts functional given by \eqref{eq:inverse_potts} for $p=2.$
We decompose $x \in \R^n$ according to  $x = x_0 + \overline{x}$, with $x_0 \in V, \overline{x} \in \R e$.
Applying this decomposition to \eqref{eq:inverse_potts} yields
\begin{align*}
   P_\gamma (x) &= \gamma \|\nabla (x_0 + \overline{x})\|_0 + \|A x_0 + A \overline{x} -b\|_2^2\\
   &= \gamma \|\nabla x_0 \|_0 + \|A x_0 + A \overline{x} -b\|_2^2.
\end{align*}
We write $\overline{x} = \mu e$ to obtain
\begin{equation} \label{eq:sparse_equiv_x}
P_\gamma(x) =    \gamma \|\nabla x_0 \|_0 + \|A  x_0 + \mu A e -b\|_2^2.
\end{equation}
Let us fix $x_0$ for the moment and let us look for $\mu = \mu(x_0)$ which minimizes
the function $\mu \to$ $P_\gamma(x_0+\mu e).$
Since $\|\nabla (x_0 + \mu e) \|_0$ $= \|\nabla (x_0 + \mu' e) \|_0$ for all $\mu,\mu'$
we have to minimize (w.r.t.\ $\mu$)
\begin{equation}
 \sum_{i=1}^m \left( \mu \sum_{j=1}^n A_{ij} + \sum_{j=1}^n A_{ij} x_{0,j} - b_i \right)^2 \to \min.
\end{equation}
The corresponding minimizer $\mu(x_0)$ can be computed explicitly (e.g., by derivating). It is given by
\begin{equation} \label{eq:mu_of_x0}
    \mu(x_0) = \frac{\sum_{i=1}^m \widetilde{A_i} b_i - \sum_{i=1}^m \widetilde{A_i} \sum_{j=1}^n A_{ij} x_{0,j}}
    { \sum_{i=1}^m \widetilde{A_i}^2},
\end{equation}
where $\widetilde{A_i}$ is the sum of the $i^{th}$ row of the matrix $A$ given by \eqref{eq:Atilde}. In particular, $\mu(x_0)$
depends affine linearly  on $x_0$, i.e. is of the form $d-Ex_0$ where $d$ is a constant and $E$ is a row vector of length $n$, both not depending on $x_0$. Plugging the expression \eqref{eq:mu_of_x0} for $\mu(x_0)$ into \eqref{eq:sparse_equiv_x}, we obtain a minimization problem in $x_0.$ It is given by
\begin{equation} \label{eq:inverse_Potts_intermediate}
    \gamma \|\nabla x_0\|_0 + \|A^\prime x_0 - b^\prime \|_2^2    \to \min,
\end{equation}
where $A^\prime$ is the matrix given by
\begin{equation} \label{eq:Aprime}
   A^\prime_{kj} = A_{kj}  - \frac{\widetilde{A_k}  \sum_{i=1}^m \widetilde{A_i} A_{ij}}{\sum_{i=1}^m \widetilde{A_i}^2},
\end{equation}
with
\begin{equation} \label{eq:Atilde}
   \widetilde{A_i} :=\sum_{j=1}^n A_{ij}
\end{equation}
and $b^\prime$ is the vector given by
\begin{equation} \label{eq:bprime}
   b_k^\prime = b_k  - \frac{ \widetilde{A_k} \sum_{i=1}^m \widetilde{A_i} b_i}{\sum_{i=1}^m \widetilde{A_i}^2}.
\end{equation}

After these preparations we show the theorem; we consider a minimizer $u^*$  of the sparsity problem \eqref{eq:sparse2}
w.r.t. the matrix $B=A'\nabla$ and data $b'$.
The crucial point is that $\nabla$ is an isomorphism from $V$ onto $\R^{n-1}$ which implies the
equivalence
\begin{equation*} 
   u^* \text{ minimizes \eqref{eq:sparse2} }  \Leftrightarrow
   x^*_0 = \nabla^+u^* \text{ minimizes \eqref{eq:inverse_Potts_intermediate} }.
\end{equation*}
Applying this equivalence, $x^*_0$ $= \nabla^+u^*$ is a minimizer of \eqref{eq:inverse_Potts_intermediate},
and, using \eqref{eq:mu_of_x0}, the vector $x^\ast$ $= x_0^* + \mu(x_0^*)$ is a minimizer of the original
Potts problem \eqref{eq:inverse_potts} for $A,b.$
\end{proof}

Using the relation between inverse Potts  and sparsity problems
we are now able to show the complexity statement \autoref{thm:Potts_NP_hard}
which asserts NP-hardness of the inverse Potts problem.
\begin{proof}[Proof of \autoref{thm:Potts_NP_hard}]
  The sparsity problem \eqref{eq:sparseLag} is NP-hard by \cite[Theorem 3]{chen2011complexity} ($p \geq 1,$ $\gamma >0.$)
  According to \autoref{thm:Pottsforsparse} each instance of the sparsity problem \eqref{eq:sparseLag}
   defines an instance  of the inverse Potts problem
  \eqref{eq:inverse_potts}. In particular, for any NP-hard instance of the sparsity problem (with matrix $A$ and data $b$)
  there is a corresponding inverse Potts problem (with matrix $A\nabla$ and data $b$.) The transformation of the
  functionals and the transformation of the corresponding minimizers given by \autoref{thm:Pottsforsparse} can obviously be done in
  polynomial time. 
  Therefore the Potts problem is NP-hard.
\end{proof}

\subsection{Convergence}\label{ssec:ProofConv}

In our presentation we have assumed that the sequence $\mu_k$ is a geometric progression.
What we actually need is that $\mu_k$ is a non-decreasing sequence fulfilling
\begin{equation} \label{eq:muk_summability}
\sum_{k} \frac{1}{\sqrt{\mu_k}}< \infty
\end{equation}
which is obviously satisfied for geometric progressions.
So we show \autoref{thm:convergence} assuming \eqref{eq:muk_summability}
instead.

\begin{proof}[Proof of \autoref{thm:convergence}]
We consider the Potts ADMM iteration for $u^k,v^k$ and $\lambda^k$ given by \eqref{eq:ADMM}.
We show that
\begin{equation}\label{eq:whatWSh}
(u^k,v^k) \to (u^*,v^*) \text{ with } u^* = v^*,  \quad\text{ and }\quad \tfrac{\lambda^k}{\mu_k} \to 0,
\end{equation}
which is a qualitative version of the assertion of the theorem.

We denote the functional occurring in the first line of \eqref{eq:ADMM} by $F_k,$ i.e.,
$$
   F_k(u) =  \gamma \|\nabla u\|_0  + \frac{\mu_k}{2} \|u - \left(v^k-\frac{\lambda^k}{\mu_k}\right)\|_2^2.
$$
Using this notation, the first line of \eqref{eq:ADMM} reads $u^{k+1} \in \argmin_u F_k(u)$.
In order to estimate $\|u^{k+1} - (v^k - \frac{\lambda^k}{\mu_k})\|_2$ we observe that
$
 F_k(u^{k+1}) \leq F_k\left(v^k-\frac{\lambda^k}{\mu_k}\right)
$
which is a consequence of the minimality of $u^{k+1}.$
Using the definition of $F_k$ yields
\begin{align*}
   \gamma \|\nabla u^{k+1}\|_0  +& \frac{\mu_k}{2} \|u^{k+1} - \left(v^k - \frac{\lambda^k}{\mu_k} \right)\|_2^2 \\
   \leq & \gamma \|\nabla \left(v^k - \frac{\lambda^k}{\mu_k}\right)\|_0 \leq \gamma n,
\end{align*}
where $n$ is the length of $v_k$.
Since the first summand on the left hand side is non-negative we get that
\begin{equation}\label{eq:est1}
\|u^{k+1} - \left(v^k - \frac{\lambda^k}{\mu_k} \right)\|_2^2 \leq \frac{\gamma n}{\mu_k}.
\end{equation}
In particular,
\begin{equation} \label{eq:lim1}
   \lim_{k \to \infty} u^{k+1} - \left(v^k - \frac{\lambda^k}{\mu_k}\right) =0.
\end{equation}

Now we draw our attention to the second line of \eqref{eq:ADMM}. We denote the corresponding functional by
$$
   G_k(v) = \|Av-b\|_p^p + \frac{\mu_k}{2}  \|v - \left(u^{k+1} + \frac{\lambda^k}{\mu_k}\right)\|_2^2.
$$
The minimality of $v^{k+1}$ implies
$G_k(v^{k+1})$ $\leq G_k\left(u^{k+1} + \frac{\lambda^k}{\mu_k}\right).$
We apply the definition of $G_k$ and estimate
\begin{align}
    \|Av^{k+1} - b\|_p^p &+ \frac{\mu_k}{2} \|v^{k+1} - \left(u^{k+1} + \frac{\lambda^k}{\mu_k}\right)\|_2^2 \notag \\ \leq & \|A\left(u^{k+1} + \frac{\lambda^k}{\mu_k}\right) - b\|_p^p \notag\\
    \leq &\|A (u^{k+1} + \frac{\lambda^k}{\mu_k} - v^k) + A v^k -b\|_p^p  \notag\\
    \leq & \left( \|A\| \|u^{k+1} + \frac{\lambda^k}{\mu_k}-v^k\|_2 +\|Av^k - b\|_p \right)^p.  \label{eq:l3}
\end{align}
Here $\|A\|$ is the norm of $A$ viewed as an operator from $\ell^2$ to $\ell^p$.
We combine the inequalities \eqref{eq:l3} and $\eqref{eq:est1}$ in order to obtain that
$$
    \|Av^{k+1} - b\|_p \leq \frac{\|A\| \gamma n}{ \mu_{k}} + \|Av^{k} - b\|_p.
    $$
Solving this recursion yields
$$
    \|Av^{k+1} - b\|_p \leq  \|A\| \gamma n \sum_{j=1}^k \frac{1}{\mu_j} + \|A v^0 -b\|_p,
$$
which shows that the sequence $(\|Av^{k+1} - b\|_p)_{k \in \N}$ is bounded. Together with \eqref{eq:l3} this implies
$$\tfrac{\mu_k}{2} \|v^{k+1} - \left(u^{k+1} +\frac{\lambda^k}{\mu_k} \right) \|_2^2 \leq
(\|A\| \|u^{k+1} + \frac{\lambda^k}{\mu_k}-v^k\|_2 + C)^p,
$$
where $C$ is a positive constant independent of $k$. Using \eqref{eq:lim1} we get that
\begin{equation} \label{eq:est2}
\mu_k \|v^{k+1} - \left(u^{k+1} +\frac{\lambda^k}{\mu_k}\right) \|_2^2 \mbox{ is bounded}.
\end{equation}

We show the convergence of the sequence $v^k$ by showing that it is a Cauchy sequence. To this end
we estimate
$$
   \|v^{k+1} - v^k\| \leq \|v^{k+1} - u^{k+1} - \frac{\lambda^k}{\mu_k}\| + \|u^{k+1} + \frac{\lambda^k}{\mu_k} - v^k\|.
$$
Now we apply \eqref{eq:est1} and \eqref{eq:est2} which yield
$$
    \|v^{k+1}-v^k\| \leq \frac{C}{\sqrt{\mu_k}}
$$
for some constant $C >0$ which is independent of $k$.
Assumption \eqref{eq:muk_summability} on $\mu_k$ guarantees that $v^k$ is a Cauchy sequence and hence that $v^k$ converges to some $v^*$.

We use the third line of \eqref{eq:ADMM} to obtain the equality
\begin{equation}\label{eq:thirdLineConsequ}
   \frac{\lambda^{k+1}}{\mu_{k+1}} = \frac{\mu_{k}}{\mu_{k+1}} \left( (\frac{\lambda^k}{\mu_k} +  u^{k+1} -v ^k)+(v^k-v^{k+1}) \right).
\end{equation}
By \eqref{eq:lim1} and \eqref{eq:est2} each term in parenthesis converges to $0.$ Since $\mu_k$ is non-decreasing, we have that
$\mu_k/\mu_{k+1} \leq 1$ and, therefore, \eqref{eq:thirdLineConsequ} implies that
$$
   \lim_{k \to \infty} \frac{\lambda^k}{\mu_k} = 0,   \quad \text{and} \quad \lim_{k \to \infty} \frac{\lambda^{k+1}}{\mu_k} = 0 .
$$
We rewrite the third line of \eqref{eq:ADMM} as $u^{k+1}- v^{k+1}$ $= (\lambda^{k+1} - \lambda^{k})/\mu_k$  
to obtain the inequality
$$
   \| u^{k+1} - v^{k+1} \| \leq \tfrac{\|\lambda^{k+1}\|} {\mu_k} + \tfrac{\|\lambda^{k}\|} {\mu_k} 
   \to 0.   
$$
This means that  $u^k - v^k \to 0$ and, since $v_k$ converges,
also $u_k$ converges and the corresponding limit $u^*$ equals $v^*$. This shows \eqref{eq:whatWSh} and completes the proof.
\end{proof}

\section{Conclusion and outlook}

We have shown that the inverse Potts problem has a minimizer in the discrete setting
but that the time continuous counterpart does not have minimizers in general.
We further have shown that the computation of minimizers is an NP-hard problem.
Having accepted that the computation of exact solutions are unfeasible, we have proposed a new approach to the inverse Potts problem based on the alternating direction method of multipliers. 
In our experiments we have compared the iPotts-ADMM algorithm with total variation minimization for jump-sparse reconstruction. We have observed that our method 
often performs better than but at least as well as TV minimization.
We further have shown that the sparsity problem can be reduced to an inverse Potts problem for $p \geq 1.$
The experiments indicate that the iPotts-based approach to the sparsity problem
performs as least as well as the state-of-the-art algorithms in presence of Gaussian noise and significantly better in presence of impulsive noise.

Future research aims at faster algorithms for the multivariate inverse Potts problem and at Potts problems with manifold valued data.

\appendices

\section{}\label{app:noise}
We consider Gaussian, Laplacian, and impulsive noise.
The first two types of noise are additive. Thus the measurement is given by
\[
	b = A \overline{x} + \eta_\sigma,
\]
where $\eta_\sigma$ is a $m$-dimensional vector of i.i.d. random variables
of standard deviation $\sigma.$
In case of Gaussian noise, the probability density function is given by
\begin{align*}
	p(x) =  \frac{1}{\sigma \sqrt{2\pi} } e^{- \frac{x^2}{2\sigma^2}}.
\end{align*}
In the case of Laplacian noise, the density is defined by
\begin{align*}
	p(x) = \frac{1}{\sigma\sqrt{2}} e^{-\frac{\sqrt{2}}{\sigma}|x|}.
\end{align*}
In the case of impulsive noise, 
we randomly choose a prescribed percentage of indices $I$ between $1$ and $n$
and set each data point belonging to that index set  to a random number, i.e.
\begin{align*}
	b_i = 
	\begin{cases}
			(A \overline{x})_i, &\text{if }i \notin I, \\
			\xi, &\text{else.}
	\end{cases}
\end{align*}
Here, $\xi$ 
is a random variable which is uniformly distributed 
in the interval $[0,1]$ for the jump-sparsity experiments
and in the interval $[-\frac12, \frac12]$ for the sparsity experiments.

\bibliographystyle{IEEEtran}
\bibliography{pottsDeconvolution}

\end{document}